\documentclass[12pt]{amsart}

\usepackage{amssymb}
\usepackage{amsmath}
\usepackage{amsfonts}
\usepackage{graphicx}
\usepackage{color}
\usepackage[onehalfspacing]{setspace}
\usepackage{hyperref}
\usepackage{natbib}
\usepackage{charter}

\makeatletter
\def\section{\@startsection{section}{1}
	\z@{1.0\linespacing\@plus\linespacing}{.8\linespacing}{\Large\centering}}

\def\subsection{\@startsection{subsection}{2}
	\z@{.8\linespacing\@plus.7\linespacing}{.7\linespacing}{\large}}

\def\subsubsection{\@startsection{subsubsection}{3}
	\z@{.5\linespacing\@plus.7\linespacing}{-.5em}{\normalfont\bfseries\centering}}
\makeatother

\numberwithin{equation}{section}

\newtheorem{theorem}{Theorem}[section]
\newtheorem{lemma}{Lemma}[section]
\newtheorem{corollary}{Corollary}[section]

\theoremstyle{definition}
\newtheorem{definition}{Definition}[section]

\theoremstyle{definition}
\newtheorem{assumption}{Assumption}[section]

\theoremstyle{definition}

\title{}
\begin{document}
	\vspace*{5ex minus 1ex}
	\begin{center}
		\Large \textsc{Stable Limit Theorems for Empirical Processes under Conditional Neighborhood Dependence}
		\bigskip
	\end{center}
	
	\date{%
		\today%
	}

	\vspace*{5ex minus 1ex}
	\begin{center}
		Ji Hyung Lee and Kyungchul Song
		
		\textit{University of Illinois and University of British Columbia}
		
		\bigskip
	\end{center}
	
	\begin{abstract}
		{\footnotesize This paper introduces a new concept of stochastic dependence among many random variables which we call conditional neighborhood dependence (CND). Suppose that there are a set of random variables and a set of sigma algebras where both sets are indexed by the same set endowed with a neighborhood system. When the set of random variables satisfies CND, any two non-adjacent sets of random variables are conditionally independent given sigma algebras having indices in one of the two sets' neighborhood. Random variables with CND include those with conditional dependency graphs and a class of Markov random fields with a global Markov property. The CND property is useful for modeling cross-sectional dependence governed by a complex, large network. This paper provides two main results. The first result is a stable central limit theorem for a sum of random variables with CND. The second result is a Donsker-type result of stable convergence of empirical processes indexed by a class of functions satisfying a certain bracketing entropy condition when the random variables satisfy CND.}\bigskip
		
        {\footnotesize \ }

\medskip

        
         \end{abstract}
	
    	\thanks{We thank Denis Kojevnikov for his valuable comments. We also thank a referee for valuable suggestions for the improvement of the paper. All errors are ours. Corresponding address: Kyungchul Song, Vancouver School of Economics, University of British Columbia, 6000 Iona Drive, Vancouver, BC, Canada, V6T 1L4. Email address: kysong@mail.ubc.ca. Co-author address: Ji Hyung Lee, Department of Economics, University of Illinois,
    		1407 W. Gregory Dr., 214 David Kinley Hall, Urbana, IL 61801. Email address: jihyung@illinois.edu}
    	\maketitle

\section{Introduction}

Empirical processes indexed by a class of functions arise in many applications, in particular in developing asymptotic inference for nonparametric or semiparametric models and for goodness-of-fit specification tests. (See, e.g., \cite{Andrews:94:Handbook} and Chapter 3 of \cite{vanderVaart/Wellner:96:WeakConvg} for a review of applications of empirical process theory in statistics and econometrics.) While a predominant body of the literature on empirical process theory focuses on independent observations or time series observations, there is relatively little research on empirical processes with spatial or cross-sectional dependence. This paper aims to contribute to the literature by providing limit theorems for empirical processes which consist of random variables with a flexible, complex (cross-sectional) dependence structure.

In this paper, we introduce a new notion of stochastic dependence among a set of random variables. Suppose that we are given a set of random variables $\{Y_i\}_{i \in N_n}$ indexed by a set $N_n$, where the set $N_n$ is endowed with a neighborhood system so that each $i \in N_n$ is associated with a subset $\nu_n(i) \subset N_n\backslash \{i\}$ called the \textit{neighborhood} of $i$. In this paper, we call the map $\nu_n: N_n \rightarrow 2^{N_n}$ a \textit{neighborhood system}. Given a neighborhood system $\nu_n$ and a set of $\sigma$-fields $\mathcal{M} \equiv (\mathcal{M}_i)_{i \in N_n}$, we say that $\{Y_i\}_{i \in N_n}$ is \textit{conditionally neighborhood dependent}(CND) with respect to $(\nu_n,\mathcal{M})$ if for any two non-adjacent subsets $A$ and $B$ of $N_n$, $((Y_i)_{i \in A}, (\mathcal{M}_i)_{i \in A})$ and $((Y_i)_{i \in B}, (\mathcal{M}_i)_{i \in B})$ are conditionally independent given $(\mathcal{M}_i)_{i \in \nu_n(A)}$, where $\nu_n(A)$ is the union of the neighborhoods of $i \in A$ with the set $A$ itself removed.

Our CND property is a generalization of both dependency graphs and Markov random fields with a global Markov property (\cite{Lauritzen/Dawid/Larsen/Leimer:1990:Networks}.) Dependency graphs were introduced by \cite{Stein:72:BerkeleySymp} in his study of normal approximation. (See \cite{Chen/Shao:04:AP} and \cite{Rinott/Rotar:96:JMA} for a general local dependence notion that is different from ours.) A set of random variables have a graph as a dependency graph, if two sets of random variables are allowed to be dependent only when the two sets are adjacent in the graph. This dependence can be viewed as restrictive in many applications, as it requires that any random variables be independent even if their indices are indirectly connected in the graph. In contrast, CND random variables are allowed to be dependent even if they are not adjacent in a graph. The CND property captures the notion that ``any two random variables are independent once we condition on the source of their joint dependence". In this sense, the CND property is closely related to a Markov property in the literature of random fields. However, in contrast to the Markov property, the CND property does not require that the $\sigma$-fields $\mathcal{M}_i$ be generated by $Y_i$ itself.

This paper provides two main results. The first main result is a Berry-Esseen bound for a sum of CND random variables. Our bound is comparable to Berry-Esseen bounds established for a sum of random variables with a dependency graph in some generic situations in the literature. (\cite{Baldi/Rinott:89:AP}, \cite{Chen/Shao:04:AP}, and \cite{Penrose:03:RGP} to name but a few.) This latter literature typically uses Stein's method to establish the bound, but to the best of our knowledge, the existing proofs using Stein's method for dependency graphs do not seem immediately extendable to a sum of CND random variables, due to a more flexible form of conditioning $\sigma$-fields involved in the CND property. In this paper, we use a traditional characteristic function-based method to derive a Berry-Esseen bound.

A typical form of a Berry-Esseen bound in these set-ups, including ours, involves the maximum degree of the neighborhood system, so that when the maximum degree is high, the bound is of little use. However, in many social networks observed, removing a small number of high-degree vertices tends to reduce the maximum degree of the neighborhood system substantially. Exploiting this insight, we provide a general version of a Berry-Esseen bound which uses conditioning on random variables associated with high degrees.

The second main result in this paper is a stable limit theorem for an empirical process indexed by a class of functions, where the empirical process is constituted by CND random variables. Stable convergence is a stronger notion of convergence than weak convergence, and is useful for asymptotic theory of statistics whose normalizing sequence has a random limit.

To obtain a stable limit theorem, we first extend the exponential inequality of \cite{Janson:04:RSA} for dependency graphs to our set-up of CND random variables, and using this, we obtain a maximal inequality for an empirical process with a bracketing-entropy type bound. This maximal inequality is useful for various purposes, especially when one needs to obtain limit theorems that are uniform over a given class of functions indexing the empirical process. Using this maximal inequality, we establish the asymptotic equicontinuity of the empirical process which, in combination with the central limit theorem that comes from our previously established Berry-Esseen bound, gives a stable limit theorem. This enables stable convergence of the empirical process to a mixture Gaussian process.

As it turns out, our stable limit theorem for an empirical process requires that the maximum degree of the neighborhood system be bounded. However, in many real-life networks, the maximum degree can be substantial, especially when the networks behave like a preferential attachment network of Barab\'{a}si-Albert.(\cite{Barabasi/Albert:99:Science}) Thus, following the same spirit of extending the Berry-Esseen bound to the case conditional on high degree vertices, we extend the stable limit theory to a set-up where it relies only on those observations with relatively low degrees by conditioning on the random variables associated with high degree vertices. This extension enables us to obtain a stable limit theorem for empirical processes when the maximum degree of the neighborhood system increases to infinity as the size of the system increases.

Stable convergence has been extensively studied in the context of martingale central limit theorems. (See, e.g. \cite{Hall/Heyde:80:MartingaleLimitTheory}.) See \cite{Hausler/Luschgy:10:PTSM} for stable limit theorems for Markov kernels and related topics. Recent studies by \cite{Kuersteiner/Prucha:2013:JOE} and \cite{Hahn/Kuersteiner/Mazzocco:16:arXiv} established a stable central limit theorem for a sum of random variables having both cross-sectional dependence and time series dependence by utilizing a martingale difference array formulation of the random variables. 

Markov-type cross-sectional dependence on a graph has received attention in the literature (\cite{Lauritzen:1996:GraphicalModels}.) In particular, the pairwise Markov property of random variables says that two non-adjacent random variables are conditionally independent given all the other variables, and is captured by a precision matrix in a high dimensional Gaussian model. (See \cite{Meinshausen/Buhlmann:08:AS} and \cite{Cai/Liu/Zhou:16:AS} for references.) This paper's CND property is stronger than the pairwise Markov property when the conditioning $\sigma$-fields, $\mathcal{M}_i$'s, are those that are generated by the random variables. However, the CND property encompasses the case where the latter condition does not hold, and thus includes dependency graphs as a special case, unlike Markov-type dependence mentioned before.

\cite{Wu:05:PNAS} introduced a dependence concept that works well with nonlinear causal processes. More recently, \cite{Jirak:16:AP} established a Berry-Esseen bound with optimal rate for nonlinear causal processes with temporal ordering. \cite{Chen/Wu:16:EJS} considered a nonlinear spatial process indexed by a lattice in the Euclidean space. These models are distinct from ours. The major distinction of our approach is to model the stochastic process to be indexed by a generic graph, and model the dependence structure using conditional independence relations along the graph. Thus our approach works well with, for example, Markov random fields on an undirected graph. On the other hand, the models of this literature accommodate various temporal or spatial autoregressive processes. To the best of our knowledge, stable convergence of empirical processes indexed by a class of functions has not been studied under either dependency graphs or Markov random fields on a graph.

The remainder of the paper proceeds as follows. In Section 2, we formally introduce the notion of conditional neighborhood dependence (CND) and study its basic properties. In Section 3, we provide stable central limit theorems for a sum of CND random variables. We also present the stable convergence of an empirical process to a mixture Gaussian process. The mathematical proofs of the results are found in the appendix.

\section{Conditional Neighborhood Dependence}
\subsection{Definition}
Let $\mathcal{N}$ be an infinite countable set. For each $n=1,2,...$, let $N_n \subset \mathcal{N}$ be a finite set such that $|N_n| = n$ and let $2^{N_n}$ be the collection of all the subsets of $N_n$. We assume that $N_n$ is a proper subset of $N_{n+1}$ for each $n \ge 1$. We will call each element of $N_n$ a \textit{vertex}, and call any map $\nu_n:N_n \rightarrow 2^{N_n}$ a \textit{neighborhood system}, if for each $i \in N_n$, $i \notin \nu_n(i)$.\footnote{Equivalently, one might view the neighborhood system as a graph by identifying each neighborhood of $i$ as the neighborhood of $i$ in the graph. However, it seems more natural to think of a stochastic dependence structure among random variables in terms of neighborhoods rather than in terms of edges in a graph.} Let us define for each $A \subset N_n$,
\begin{eqnarray*}
	\bar \nu_n(A) \equiv \left(\bigcup_{i \in A} \nu_n(i) \right) \cup A \textnormal{ and } 
	\nu_n(A) \equiv \bar \nu_n(A) \setminus A.
\end{eqnarray*} 
For $\{i_1,...,i_m\} \subset N_n$, we simply write $\nu_n(i_1,...,i_m) = \nu_n(\{i_1,...,i_m\})$ and $\bar \nu_n(i_1,...,i_m) = \bar \nu_n(\{i_1,...,i_m\})$, suppressing the curly brackets. Let us call $\bar \nu_n (A)$ \textit{the $\nu_n$-closure of} $A$ and $\nu_n(A)$ \textit{the $\nu_n$-boundary of} $A$. The $\nu_n$-closure of $A$ includes the vertices in $A$ but the $\nu_n$-boundary around $A$ excludes them.

If for each $i,j \in N_n$, $i \in \nu_n(j)$ implies $j \in \nu_n(i)$, we say that neighborhood system $\nu_n$ is \textit{undirected}. If there exists a pair $i,j \in N_n$ such that $i \in \nu_n(j)$ but $j \notin \nu_n(i)$, we say that neighborhood system $\nu_n$ is \textit{directed}.

It is often useful to compare different dependence structures governed by different neighborhood systems. When we have two neighborhood systems $\nu_n$ and $\nu_n'$ such that $\nu_n(i) \subset \nu_n'(i)$ for each $i \in N_n$ and $\nu_n(j) \ne \nu_n'(j)$ for some $j \in N_n$, we say that $\nu_n$ \textit{is strictly finer than} $\nu_n'$ and $\nu_n'$ \textit{is strictly coarser than} $\nu_n$. When $\nu_n(i) \subset \nu_n'(i)$ for each $i \in N_n$, we say that $\nu_n$ \textit{is weakly finer than} $\nu_n'$ and $\nu_n'$ \textit{is weakly coarser than} $\nu_n$.

Let us introduce the notion of dependence among a triangular array of $\sigma$-fields. Let $(\Omega,\mathcal{F},P)$ be a given probability space, and let  $\{\mathcal{F}_{i}\}_{i \in N_n}$ be a given triangular array of sub-$\sigma$-fields of $\mathcal{F}$, indexed by $i \in N_n$. (Proper notation for the sub-$\sigma$-field in the triangular array should be $\mathcal{F}_{in}$, but we suppress the $n$ subscript for simplicity.) For any $A \subset N_n$, we let $\mathcal{F}_A$ be the smallest $\sigma$-field that contains $\mathcal{F}_{i},i \in A$, and $\mathcal{F}_{-A}$ be the smallest $\sigma$-field that contains all the $\sigma$-fields $\mathcal{F}_{i}$ such that $i \in N_n \backslash A$. When $A = \varnothing$, we simply take $\mathcal{F}_A$ to be the trivial $\sigma$-field. We apply this notation to other triangular arrays of $\sigma$-fields, so that if $\{\mathcal{M}_{i}\}_{i \in N_n}$ is a triangular array of sub-$\sigma$-fields of $\mathcal{F}$, we similarly define $\mathcal{M}_{A}$ and $\mathcal{M}_{-A}$ for any $A \subset N_n$. For given two $\sigma$-fields, say, $\mathcal{G}_1$ and $\mathcal{G}_2$, we write $\mathcal{G}_1 \vee \mathcal{G}_2$ to represent the smallest $\sigma$-field that contains both $\mathcal{G}_1$ and $\mathcal{G}_2$.

Given a triangular array of $\sigma$-fields, $\{\mathcal{M}_i\}_{i=1}^\infty$, let us introduce a sub $\sigma$-field $\mathcal{G}$ defined by
\begin{eqnarray}
\label{sigma G}
	\mathcal{G} \equiv \bigcap_{n \ge 1} \bigcap_{i \in N_n} \mathcal{M}_i.
\end{eqnarray}
In many applications, $\mathcal{G}$ is used to accommodate random variables with a common shock. For example, suppose that each $\mathcal{M}_i$ is generated by a random vector $(\varepsilon_i,U)$ from a set of random variables $\{\varepsilon_i\}_{i \in N_n}$ that are conditionally independent given a common random variable $U$. Then we can take $\mathcal{G}$ to be the $\sigma$-field generated by $U$. We will discuss examples of CND random vectors in a later section, after we study their properties.

Let us introduce the notion of dependence of an array of $\sigma$-fields that is of central focus in this paper.

\begin{definition}	
	\noindent (i) Given neighborhood system $\nu_n$ on $N_n$ and an array of $\sigma$-fields, $\mathcal{M} \equiv \{\mathcal{M}_{i}\}_{i \in N_n}$, we say that $\sigma$-fields  $\{\mathcal{F}_{i}\}_{i \in N_n'}$ for a given subset $N_n' \subset N_n$ are  \textit{conditionally neighborhood dependent (CND) with respect to} $(\nu_n,\mathcal{M})$, if for any $A,B \subset N_n'$ such that $A \subset N_n' \setminus \bar \nu_n(B)$ and $B \subset N_n' \setminus \bar \nu_n(A)$, $\mathcal{F}_A \vee \mathcal{M}_A$ and $\mathcal{F}_B \vee \mathcal{M}_B$ are conditionally independent given $\mathcal{M}_{\nu_n(A)}$.
    \medskip
    	
	\noindent (ii) If $\sigma$-fields generated by random vectors in $\{Y_{i}\}_{i \in N_n'}$ for a subset $N_n' \subset N_n$ are CND with respect to  $(\nu_n,\mathcal{M})$, we simply say that random vectors in $\{Y_{i}\}_{i \in N_n}$ are \textit{CND with respect to} $(\nu_n,\mathcal{M})$.
\end{definition}

Conditional neighborhood dependence specifies only how conditional independence arises, not how conditional dependence arises. Conditional neighborhood dependence does not specify independence or dependence between $Y_i$ and $Y_j$ if $j \in \nu_n(i)$ or $i \in \nu_n(j)$. Furthermore, conditional neighborhood dependence can accommodate the situation where the neighborhoods in the system $\nu_n$ are generated by some random graph on $N_n$, as long as the random graph is $\mathcal{G}$-measurable. In such a situation, the results of this paper continue to hold with only minor modifications that take care of the randomness of $\nu_n$.

\begin{figure}[t]
	\centering
	
	\label{fig:CND}
	\vspace{0.5cm}
	\includegraphics[scale=0.6]{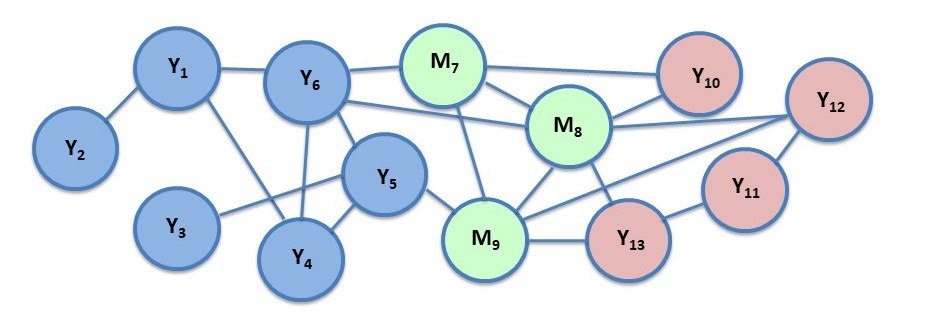}
	\caption{Conditional Neighborhood Dependence}
	\medskip
	
	\parbox{6.2in}{\footnotesize
		Notes: Suppose we are given random variables $(Y_1,...,Y_{13})$ and $(M_1,...,M_{13})$, where $\mathcal{M}_i = \sigma(M_i)$. The figure depicts the neighborhood system on $N_n = \{1,2,...,13\}$. If $Y_1,...,Y_{13}$ are CND with respect to the neighborhood system as illustrated above and $(\mathcal{M}_1,...,\mathcal{M}_{13})$, it implies that $Y_1,...,Y_6$ and $Y_{10},...,Y_{13}$ are conditionally independent given $\mathcal{M}_7,\mathcal{M}_8,\mathcal{M}_9$.\bigskip \bigskip \bigskip} 
\end{figure}

\subsection{Monotonicity and Invariance to Conditioning}

In general, a conditional independence property is not monotone in conditioning $\sigma$-fields. In other words, when $X$ and $Y$ are conditionally independent given a $\sigma$-field $\mathcal{F}'$, this does not imply conditional independence between $X$ and $Y$ given any sub-$\sigma$-field of $\mathcal{F}'$ or given any $\sigma$-field that contains $\mathcal{F}'$. However, CND partially obeys monotonicity in neighborhood systems. More specifically, the CND property with a finer neighborhood system implies more conditional independence restrictions than what the CND property with a coarser neighborhood system implies. We introduce a lemma which makes precise this monotonicity property of CND. Suppose that we are given two neighborhood systems $\nu_n$ and $\nu_n'$ where $\nu_n$ is weakly finer than $\nu_n'$. The following lemma shows that the CND property of a triangular array with respect to a given neighborhood system $\nu_n'$ carries over to that with respect to $\nu_n$ \textit{if} $\nu_n$ is undirected.

\begin{lemma} 
	\label{lemm: monotonicity}
Given two neighborhood systems $\nu_n$ and $\nu_n'$, suppose that $\nu_n$ is weakly finer than $\nu_n'$ and that a triangular array of $\sigma$-fields $\{\mathcal{F}_i\}_{i \in N_n}$ is CND with respect to $(\nu_n,\mathcal{M})$ for some $\sigma$-fields $\mathcal{M}\equiv\{\mathcal{M}_{i}\}_{i \in N_n}$. Suppose further that $\nu_n$ is undirected. Then $\{\mathcal{F}_i\}_{i \in N_n}$ is CND with respect to $(\nu_n',\mathcal{M})$.
\end{lemma}

\noindent \textbf{Proof: } Take any $A,B \subset N_n$ such that $A \subset N_n \setminus \bar \nu_n'(B)$ and $B \subset N_n \setminus \bar \nu_n'(A)$, so that $A \subset N_n \setminus \bar \nu_n(B)$ and $B \subset N_n \setminus \bar \nu_n(A)$ as well. Let $D \equiv \nu_n'(A) \setminus \bar \nu_n(A)$, which, by the undirectedness of $\nu_n$, implies that $A \subset N_n \setminus \bar \nu_n(D)$. Therefore, $A \subset N_n \setminus \bar \nu_n(B \cup D)$ and $B \cup D \subset N_n \setminus \bar \nu_n(A)$. By the CND property, $\mathcal{F}_{A} \vee \mathcal{M}_{A}$ and $\mathcal{F}_{B \cup D} \vee \mathcal{M}_{B \cup D}$ are conditionally independent given $\mathcal{M}_{\nu_n(A)}$. By Lemma 4.2 of \cite{Dawid:79:JRSS}, this implies that $\mathcal{F}_{A} \vee \mathcal{M}_A$ and $\mathcal{F}_{B} \vee \mathcal{M}_{B}$ are conditionally independent given $\mathcal{M}_{\nu_n(A) \cup D}$. This proves the lemma, because 
\begin{eqnarray*}
	\nu_n(A) \cup D = \nu_n'(A),
\end{eqnarray*}
which follows due to $\nu_n(A) \subset \nu_n'(A)$, $\nu_n(A) \subset \bar \nu_n(A)$, and $\bar \nu_n(A) \setminus \nu_n(A) = A$. $\blacksquare$
\medskip

In the above lemma, the requirement that $\nu_n$ be undirected cannot be eliminated. To see this, consider the following counterexample of the lemma when $\nu_n$ is taken to be a directed neighborhood system.
\medskip

\noindent \textbf{Example 1:} Let us take $N_n=\{1,2,3,4\}$ with $n=4$. Let $(\varepsilon_i)_{i \in N_n}$ be i.i.d. standard normal random variables. Let us take
\begin{eqnarray*}
	Y_1 = \varepsilon_1, Y_2 = \varepsilon_2, Y_3 = \varepsilon_1+\varepsilon_3+\varepsilon_4, \textnormal{ and } Y_4 = \varepsilon_4. 
\end{eqnarray*}
Let a neighborhood system $\nu_n$ be given as $\nu_n(i) = \varnothing$ for all $i=1,2,4,$ and $\nu_n(3) = \{1,4\}$. Hence $\nu_n$ is directed. Then we take $\mathcal{M} = \mathcal{F} = (\mathcal{F}_i)_{i \in N_n}$, where each $\mathcal{F}_i$ denotes the $\sigma$-field generated by $Y_i$. Then it is not hard to see that $(Y_i)_{i \in N_n}$ is CND with respect to $(\nu_n,\mathcal{M})$.

Now let us introduce another neighborhood system $\nu_n'$ that is weakly coarser than $\nu_n$. Let $\nu_n'(1) = \{3\}$, $\nu_n'(2) = \{1\}$, $\nu_n'(3) = \{1,4\}$, and $\nu_n'(4) = \{3\}$. Then we take $A=\{1\}$ and $B=\{4\}$ so that $A \subset N_n \setminus \bar \nu_n'(B)$ and $B \subset N_n \setminus \bar \nu_n'(A)$. Note that 
$\nu_n'(A) = \nu_n'(1) = \{3\}$. Certainly $Y_1$ and $Y_4$ are not conditionally independent given $Y_3$, because $Y_3$ involves both $\varepsilon_1$ and $\varepsilon_4$. $\blacksquare$
 \medskip

As we shall see later, using further conditioning in the CND property, one may obtain a better normal approximation for a sum of CND random variables in some situations. Here we give a preliminary result which addresses the question of whether the CND property still holds after we increase the $\sigma$-fields in $\mathcal{M}$ in a certain way. More precisely, suppose that $Y \equiv \{Y_i\}_{i \in N_n}$ is CND with respect to $(\nu_n,\mathcal{M})$. Then let us choose $N_n^* \subset N_n$, and define a neighborhood system $\nu_n^*$ on $N_n$ such that for all $A \subset N_n$,
\begin{eqnarray}
\label{nu star}
	\nu_n^*(A) \equiv \nu_n(A) \cap N_n^*.
\end{eqnarray}
Let $\mathcal{M}^* \equiv \{\mathcal{M}_i^*\}_{i \in N_n}$, where
\begin{eqnarray}
\label{M star}
	\mathcal{M}_i^* \equiv \mathcal{M}_i \vee \mathcal{M}_{N_n \backslash N_n^*}.
\end{eqnarray}
Then we obtain the following result.
 
\begin{lemma}
	\label{lemm: invar}
	Suppose that $Y\equiv \{Y_i\}_{i \in N_n}$ is CND with respect to $(\nu_n,\mathcal{M})$. Fix $N_n^* \subset N_n$ and define $\nu_n^*$ as in (\ref{nu star}) and $\mathcal{M}^* \equiv \{\mathcal{M}_i^*\}_{i \in N_n}$ as in (\ref{M star}). Then $Y^* \equiv \{Y_i\}_{i \in N_n^*}$ is CND with respect to $(\nu_n^*,\mathcal{M}^*)$.
\end{lemma}

\noindent \textbf{Proof: } Take $A,B \subset N_n^*$ such that $A \subset N_n^* \backslash \bar \nu_n(B)$ and $B \subset N_n^* \backslash \bar \nu_n(A)$. Note that for any $A' \subset N_n^*$, $N_n^* \backslash \bar \nu_n(A') = N_n^* \backslash \bar \nu_n^*(A')$. Let
\begin{eqnarray*}
	B^* \equiv B \cup ((N_n \backslash N_n^*) \backslash \bar \nu_n(A)).
\end{eqnarray*}
By the CND of $Y$, we have that $\sigma(Y_A)\vee \mathcal{M}_A$ and $\sigma(Y_{B^*})\vee \mathcal{M}_{B^*}$ are conditionally independent given $\mathcal{M}_{\nu_n(A)}$. Hence $\sigma(Y_A)\vee \mathcal{M}_A$ and $\sigma(Y_B)\vee \mathcal{M}_{B^*}$ are conditionally independent given $\mathcal{M}_{\nu_n(A)}$. From (\ref{nu star}), we have
\begin{eqnarray*}
	\nu_n^*(A)^c \cap \nu_n(A) = \nu_n(A) \cap (N_n \backslash N_n^*).
\end{eqnarray*}
Hence we can write
\begin{eqnarray*}
	\nu_n(A) = \nu_n^*(A) \cup (\nu_n^*(A)^c \cap \nu_n(A))
	= \nu_n^*(A) \cup (\nu_n(A) \cap (N_n \backslash N_n^*)),
\end{eqnarray*}
so that $\sigma(Y_A)\vee \mathcal{M}_A$ and $\sigma(Y_B)\vee \mathcal{M}_{B^*}$ are conditionally independent when we condition on  $\mathcal{M}_{\nu_n^*(A) \cup (\nu_n(A) \cap (N_n \backslash N_n^*))}$. This implies that $\sigma(Y_A)\vee \mathcal{M}_A$ and $\sigma(Y_B)\vee \mathcal{M}_B$ are conditionally independent given  $\mathcal{M}_{\nu_n^*(A) \cup (N_n \backslash N_n^*)}$ by the definition of $B^*$, and by the choice of $A \subset N_n^*$. $\blacksquare$
\medskip
 
\subsection{Examples}

\subsubsection{Conditional Dependency Graphs}
Let $G_n=(N_n,E_n)$ be an undirected graph on $N_n$, where $E_n$ denotes the set of edges. Define $\nu_n(i) = \{j \in N_n: ij \in E_n\}$. Suppose that $\{Y_i\}$ has $G_n$ as a \textit{conditional dependency graph}, i.e., for each $A \subset N_n$, $Y_{A}$ and $Y_{N_n\backslash \bar \nu_n(A)}$ are conditionally independent given a $\sigma$-field $\mathcal{C}$. Dependency graphs were introduced by \cite{Stein:72:BerkeleySymp} and have received attention in the literature. (See, for example, \cite{Janson:88:AP}, \citet*{Baldi/Rinott:89:AP} and \citet*{Rinott/Rotar:96:JMA} for Berry-Esseen bounds for the sum of random variables with a dependency graph, and \cite{Janson:04:RSA} for an exponential inequality. See \citet*{Song:15:WP} for an application to permutation inference.) Then $\{Y_i\}$ is CND with respect to $(\nu_n,\mathcal{M})$, where $\mathcal{M}_{i}$ is taken to be $\mathcal{C}$ for each $i \in N_n$. $\blacksquare$

\subsubsection{Functional Local Dependence}
\label{sec: functional local dependence}
Let $G_n=(N_n,E_n)$ be a given directed graph, so that $ij \in E_n$ represents an edge \textit{from } $i$ \textit{to} $j$. The neighborhood $N_n^O(i) = \{j \in N_n: ij \in E_n\}$ is the \textit{out-neighborhood} of $i$, the set of vertices $j$ that there is an edge from vertex $i$ to vertex $j$. Similarly, $N_n^I(i) = \{j \in N_n: ji \in E_n\}$ is the \textit{in-neighborhood} of $i$. We define $\bar{N_n}^O(i) = N_n^O(i) \cup \{i\}$ and $\bar{N_n}^I(i) = N_n^I(i) \cup \{i\}$. Suppose that $Y_{i}$ is generated in the following way:
\begin{eqnarray}
\label{index_Y}
Y_{i} = \gamma_i(\varepsilon_{\bar N_n^I (i)},\eta_i),
\end{eqnarray} 
where $\varepsilon_{\bar N_n^I (i)} = (\varepsilon_j)_{j \in \bar N_n^I (i)}$ and $\varepsilon_i$'s and $\eta_i$'s are independent across $i$'s and $(\varepsilon_i)_{i \in N_n}$ and $(\eta_i)_{i \in N_n}$ are independent from each other. Here the functions, $\gamma_i$'s, are nonstochastic.

This way of modeling local dependence among $Y_i$'s through base random variables, $\varepsilon_i$'s, is useful in many contexts of applications. In particular, each outcome $Y_{i}$ may arise as a consequence of local interactions among individual variables, where locality is determined by a given graph $G_n$. (See \cite{Leung:16:WP} and \cite{Canen/Schwartz/Song:17:WP} for applications in economics.)

Let us define a new (undirected) graph $G_n'=(N_n,E_n')$ such that $ij \in E_n'$ if and only if $\bar N_n^I(i) \cap \bar N_n^I(j) \ne \varnothing$. In other words, $i$ and $j$ are adjacent in $G_n'$, if their in-neighborhoods overlap. Define $N_n'(i) = \{j \in N_n: ij \in E_n'\}$ for each $i \in N_n$. Then it is not hard to see that $\{Y_{i}\}_{i \in N_n}$ is CND with respect to $(N_n',\mathcal{M}')$, for any triangular arrays of $\sigma$-fields $\mathcal{M}_i'$, as long as each $\mathcal{M}_i'$ is
the $\sigma$-fields generated by $f_i((\varepsilon_j,\eta_j)_{j \in N_n})$ for any nonstochastic measurable map $f_i$. Note that we can take $\mathcal{M}_i'$ to be trivial $\sigma$-fields in which case the functional local dependence is essentially equivalent to the dependency graph assumption.

However, functional local dependence has much richer implications than the dependency graph assumption alone, because it generates lots of conditional independence restrictions that are not implied by the dependency graph assumption alone. Using such restrictions, we can obtain conditional neighborhood dependence as follows: for each $i \in N_n$, let $\mathcal{M}_{i}$ be the $\sigma$-field generated by $(\varepsilon_i,\eta_i)$. It is not hard to see that the triangular array $\{Y_{i}\}_{i \in N_n}$ is CND with respect to $(N_n^I,\mathcal{M})$. Note that the neighborhood system $N_n^I$ is weakly finer than $N_n'$, and hence when the graph $G_n$ is undirected, by Lemma \ref{lemm: monotonicity}, the CND with respect to $(N_n^I,\mathcal{M})$ expresses richer conditional independence restrictions than the CND with respect to $(N_n',\mathcal{M})$.

The notion of functional local dependence is related to physical dependence of \citet*{Wu:05:PNAS}. The difference is that physical dependence in \citet*{Wu:05:PNAS} is mainly intended for time series dependence, where $Y_{i}$ involves only the past and present values of $\varepsilon_{i}$'s through a common function, whereas the functional local dependence captures local dependence through a system of neighborhoods of random variables. $\blacksquare$

\subsubsection{Markov Random Fields on a Neighborhood System}
Suppose that we have a triangular array of random vectors $\{Y_i\}_{i \in N_n}$, $Y_i \in \mathbf{R}^d$, where there is an \textit{undirected} neighborhood system $\nu_n$ on $N_n$. A \textit{path} in $\nu_n$ between two vertices $i$ and $j$ is defined to be a sequence of distinct vertices $i_1, i_2,...,i_m$ such that $i_1 = i$ and $i_m = j$ and $i_t \in \nu_n(i_{t+1})$ for each $t=1,...,m-1$. Let us say that a set $S \subset N_n$ \textit{separates} sets $A,B \subset N_n$, if for any $i \in A$ and any $j \in B$, every path between a vertex in $A$ and a vertex in $B$ intersects set $S$. Let us consider the following two notions of Markov properties (see \cite{Lauritzen:1996:GraphicalModels}, p.32).
\begin{definition}
	(i) We say that $\{Y_i\}_{i \in N_n}$ satisfies the \textit{the pairwise Markov property} if for any two vertices $i,j \in N_n$ such that $i \notin \overline \nu_n(j)$ and $j \notin \overline \nu_n(i)$, $Y_i$ and $Y_j$ are conditionally independent given $Y_{N \setminus \{i,j\}}$.
	
	(ii) We say that $\{Y_i\}_{i \in N_n}$ satisfies the \textit{the local Markov property} if for any vertex $i \in N_n$ and any $B \subset N_n\backslash \bar \nu_n(i)$, $Y_i$ and $Y_B$ are conditionally independent given $Y_{\nu_n(i)}$.
	
	(iii) We say that $\{Y_i\}_{i \in N_n}$ satisfies the \textit{the global Markov property} if for any two subsets $A$ and $B$ of $N_n$ which are separated by set $S \subset N_n$, $Y_A$ and $Y_B$ are conditionally independent given $Y_S$.
\end{definition}
Suppose that $\mathcal{M}_i$ is the $\sigma$-field generated by $Y_i$ for each $i \in N_n$. It is not hard to see that if $\{Y_i\}_{i \in N_n}$ satisfies the global Markov property, it is CND with respect to $(\nu_n,\mathcal{M})$. And if $\{Y_i\}_{i \in N_n}$ is CND with respect to $(\nu_n,\mathcal{M})$, it satisfies the local Markov property. Hence our notion of CND is an intermediate concept between the local and the global Markov properties.\footnote{
\cite{Lauritzen/Dawid/Larsen/Leimer:1990:Networks} proposed Markov fields over directed acyclic graphs. They defined local Markov property and global Markov property and provided a sufficient condition under which both are equivalent.} Suppose that for each $i \in N$, $Y_i$ takes values in a finite set, say, $\mathcal{Y}_i$, and $Y = (Y_i)_{i \in N_n}$ is a discrete random vector with a positive probability mass at each point in the Cartesian product $\times_{i \in N}\mathcal{Y}_i$. Then, the pairwise Markov propety implies the global Markov property, and hence implies the CND property. (See p.119 of \cite{Koller/Friedman:09:ProbGraph}.)

Note that Markov chains are not CND in general. For example, consider the set $N_n = \{1,2,...,n\}$ which represents time and a directed graph $\nu_n$ on $N_n$ such that $\nu_n(i) = \{i-1\}$. Let $\{Y_i\}_{i \in N_n}$ be a Markov chain. Then the requirement from a CND property that $Y_i$ and $Y_{i+2}$ be conditionally independent given $Y_{i-1}$ does not follow from the Markov chain property.

\section{Stable Limit Theorems} 
\subsection{Stable Central Limit Theorems}

\subsubsection{The Basic Result}
In this section, we give a Berry-Esseen bound conditional on $\mathcal{G}$ for a triangular array of random variables that are CND with respect to $(\nu_n,\mathcal{M})$. Given a neighborhood system $\nu_n$ on $N_n$, we define
\begin{eqnarray*}
	d_{mx} \equiv \max_{i \in N_n}|\nu_n(i)| \textnormal{ and } d_{av} \equiv \frac{1}{n} \sum_{i \in N_n} |\nu_n(i)|,
\end{eqnarray*}
where $|\nu_n(i)|$ denotes the cardinality of set $\nu_n(i)$. We call $d_{mx}$ the \textit{maximum degree} and $d_{av}$ the \textit{average degree} of neighborhood system $\nu_n$. We use $d_{mx}$ and $d_{av}$ to express the conditions for the neighborhood system $\nu_n$. For $p=1,...,4$, define
\begin{eqnarray*}
	\mu_{p} \equiv \max_{i \in N_n} \left(\mathbf{E}\left[\left|\frac{X_{i}}{\sigma_n}\right|^p|\mathcal{G}\right] \right)^{1/p},
\end{eqnarray*}
where $\sigma_n^2 \equiv Var(\sum_{i \in N_n}X_{i}|\mathcal{G})$. Let
\begin{eqnarray*}
	r_n^2 \equiv \frac{1}{\sigma_n^4} \mathbf{E}\left[\left(\sum_{i \in N_n} \sum_{j \in \bar \nu_n(i)} \xi_{ij} \right)^2 |\mathcal{G}\right],
\end{eqnarray*}
where
$\xi_{ij} \equiv \mathbf{E}[X_iX_j|\mathcal{M}_{\nu_n(i,j)}]-\mathbf{E}[X_iX_j|\mathcal{G}]$.

\begin{theorem}
	\label{thm: CLT}
	Suppose that a triangular array $\{X_{i}\}_{i \in N_n}$ is CND with respect to $(\nu_n,\mathcal{M})$. Furthermore assume that $\mathbf{E}[X_{i}|\mathcal{M}_{\nu_n(i)}]=0$, a.e. for each $i \in N_n$, and let
	\begin{eqnarray}
		\label{event}
		\mathcal{A}_n \equiv\{\omega \in \Omega: n d_{mx} d_{av}\mu_3^3(\omega) \le 1 \}.
	\end{eqnarray}
	
	Then there exists an absolute constant $C>0$ such that on the event $\mathcal{A}_n$,  for each $n \ge 1$,
	\begin{eqnarray}
	\label{bd5}
	\sup_{t \in \mathbf{R}} \Delta_n(t;\mathcal{G}) \le C \left(\sqrt{n d_{mx} d_{av}\mu_3^3} - \log(n d_{mx} d_{av}\mu_3^3) \sqrt{n d_{mx}^2 d_{av} \mu_4^4 + r_n^2 } \right), \text{  a.e.,} \notag
	\end{eqnarray}
	where
	\begin{eqnarray*}
		\Delta_n(t;\mathcal{G}) &\equiv& \left|P\left\{\frac{1}{\sigma_n} \sum_{i \in N_n}  X_{i} \le t|\mathcal{G} \right\}  - \Phi(t)\right|
	\end{eqnarray*}
	and $\Phi$ is the distribution function of the standard normal distribution.
\end{theorem}

Since the conditional CDF is a.e. right continuous, $\sup_{t \in \mathbf{R}} \Delta_n(t;\mathcal{G}) = \sup_{t \in \mathbf{Q}} \Delta_n(t;\mathcal{G})$, a.e., for any countable dense subset $\mathbf{Q}$ of $\mathbf{R}$. The use of the bound requires  a good bound for $r_n$. Observe that $\mathbf{E}[\xi_{ij}|\mathcal{G}] = 0$. Therefore, when $\xi_{ij}$'s are locally dependent in a proper sense, we can expect that $r_n^2$ is at most of the same order as the term $nd_{mx}^2 d_{av} \mu_4^4$. The following corollary gives a set of conditions under which this is true.
\medskip

\begin{lemma}
	\label{CLT cond}
	Suppose that the conditions of Theorem \ref{thm: CLT} hold. Furthermore, $\mathcal{M}_i$'s are conditionally independent given $\mathcal{G}$. Then
	\begin{eqnarray*}
		r_n^2 \le 8 n d_{mx}^2 d_{av} \mu_4^4. 
	\end{eqnarray*}
\end{lemma}

Focusing on a special case satisfying an additional condition below, we can obtain an improved version of Theorem \ref{thm: CLT}.
\medskip

\noindent \textbf{Condition A:} For any $A, B \subset N_n$ such that $A \subset N_n \setminus \bar \nu_n(B)$, $\mathcal{M}_{\nu_n(A)}$ and $\mathcal{F}_{B}$ are conditionally independent given $\mathcal{G}$.
\medskip

Condition A accommodates conditional dependency graphs but excludes Markov random fields.

\begin{corollary}
	\label{CLT2}
	Suppose that a triangular array $\{X_{i}\}_{i \in N_n}$ is CND with respect to $(\nu_n,\mathcal{M})$, and that Condition A holds. Suppose further that $\mathbf{E}[X_{i}|\mathcal{G}]=0$, a.e. for each $i \in N_n$. 
	
	Then there exists an absolute constant $C>0$ such that on the event $\mathcal{A}_n$ defined in (\ref{event}), for each $n \ge 1$,
	\begin{eqnarray}
	\label{bd}
	\sup_{t \in \mathbf{R}} \Delta_n(t;\mathcal{G}) \le C\left(\sqrt{nd_{mx} d_{av}\mu_{3}^3}
	- \log(n d_{mx} d_{av}\mu_3^3)\sqrt{n d_{mx}^2 d_{av} \mu_{4}^4} \right), \text{  a.e.,} \notag
	\end{eqnarray}
	where we define $\mu_{p}$ and $\Delta_n(t;\mathcal{G})$ as in Theorem \ref{thm: CLT}.
\end{corollary}

The improvement of the result due to Condition A is two fold. First, the condition $\mathbf{E}[X_i|\mathcal{M}_{\nu_n(i,j)}] = 0$ is weakened to  $\mathbf{E}[X_i|\mathcal{G}] = 0$. Second, the bound does not involve $r_n^2$. When it is cumbersome to compute a reliable bound for $r_n^2$, the above corollary can be useful.\footnote{In the special case of conditional dependency graphs, one can follow the proof of Theorem 2.4 of \cite{Penrose:03:RGP} to obtain a slightly improved bound that does not have the logarithmic factor. It appears that this improvement is marginal in many applications. For example, the quantity $\sqrt{n d_{mx}^2 d_{av} \mu_{4}^4}$ is asymptotically dominated by $\sqrt{nd_{mx} d_{av}\mu_{3}^3}$, when $\mu_p^p = O_P(n^{-p/2})$ and $d_{mx}$ increases with $n$ slower than the rate $n^{1/2}$.}

In the case of dependency graphs, there has been much research establishing a Berry-Esseen bound. When we confine our attention to the special case of $|X_i| \le 1/\sqrt{n}$, $i \in N_n$, and $d_{mx} <D$ for all $n \ge 1$ for some constant $D>0$, our bound in Corollary \ref{CLT2} has the same $n^{-1/2}$ rate as in \cite{Baldi/Rinott:89:AP} (Corollary 2), \cite{Chen/Shao:04:AP} (Theorem 2.7) and \cite{Penrose:03:RGP} (Theorem 2.4), among others. These papers adopted Stein's method to obtain the bound. However, to the best of our knowledge, it is not straightforward to extend their results to our set-up of CND variables. The main reason is that the conditioning $\sigma$-field conditional on which two sets of random variables $X_A$ and $X_B$ are independent varies depending on the set $A$. Thus, for example, we cannot apply Equation (2.4) in \cite{Penrose:03:RGP}, p.31, in our context. In this paper, we resort to a more traditional Fourier analytic method in combination with Esseen's inequality.

The Berry-Esseen bound gives stable convergence of a sum of random vectors to a mixture normal distribution. More specifically, suppose that $\{X_i\}_{i \in N_n}$ is a triangular array of random variables such that $\mathbf{E}[X_i|\mathcal{M}_{\nu_n(i)}] = 0$, a.e., for each $i \in N_n$, and for each $t \in \mathbf{R}$,
	\begin{eqnarray*}
		\Delta_n(t;\mathcal{G}) \rightarrow_P 0, \text{ as } n \rightarrow \infty.
	\end{eqnarray*}
	Then for each $U \in \mathcal{G}$, and for each uniformly continuous and bounded map $f$ on $\mathbf{R}$, we have
	\begin{eqnarray*}
		\mathbf{E}\left[f\left(\frac{1}{\sigma_n} \sum_{i \in N_n}  X_{i} \right) 1_U \right] \rightarrow \mathbf{E}[f(\mathbb{Z}) 1_U], \textit{ as } n \rightarrow \infty,
	\end{eqnarray*}
	where $\mathbb{Z}$ is a standard normal random variable that is independent of $U$ and $1_U$ denotes the indicator of event $U$.

\subsubsection{Conditional Neighborhood Dependence Conditional on High-Degree Vertices}

Let us extend Theorem \ref{thm: CLT} by considering normal approximation conditioning on random variables associated with high degree vertices.

Let $N_n^* \subset N_n$ be a given subset, and $\nu_n^*$ the neighborhood system on $N_n$ given as in (\ref{nu star}). Let
\begin{eqnarray}
\label{max deg star}
	d_{mx}^* \equiv \max_{i \in N_n^*} |\nu_n^*(i)|, \text{ and }
	d_{av}^* \equiv \frac{1}{n^*} \sum_{ j \in N_n^*} |\nu_n^*(i)|,
\end{eqnarray}
where $n^* \equiv |N_n^*|$. Hence, $d_{mx}^*$ and $d_{av}^*$ are the maximum and average degrees of the restriction of $\nu_n$ to $N_n^*$. Moreover, define $\mathcal{M}^* = (\mathcal{M}_i^*)_{i \in N_n}$, where $\mathcal{M}_i^* = \mathcal{M}_i \vee \mathcal{M}_{N_n \setminus N_n^*}$, and write $\mathcal{M}^*_{\nu_n^*(i)} = \mathcal{M}_{\nu_n^*(i)} \vee \mathcal{M}_{N_n \backslash N_n^*}$. We have in mind choosing $N_n^*$ so that the set $N_n \backslash N_n^*$ consists only of high-degree vertices in the neighborhood system $\nu_n$. Then by Lemma \ref{lemm: invar}, if $(X_i)_{i \in N_n}$ is CND with respect to $(\nu_n,\mathcal{M})$, then $(X_i)_{i \in N_n^*}$ is CND with respect to $(\nu_n^*,\mathcal{M}^*)$. The main idea is that if the difference between the two sums
\begin{eqnarray}
\label{two sums}
	\sum_{i \in N_n^*} (X_i - \mathbf{E}[X_i|\mathcal{M}^*_{\nu_n^*(i)}]) \text{ and } \sum_{i \in N_n} (X_i - \mathbf{E}[X_i|\mathcal{M}_{\nu_n(i)}])
\end{eqnarray}
is asymptotically negligible, we can use the Berry-Esseen bound for the first sum using Theorem \ref{thm: CLT} and deal with the remainder term that comes from the difference.

Now let us present an extended version of Theorem \ref{thm: CLT}. Let
\begin{eqnarray*}
	\mu^*_p \equiv \max_{i \in N_n^*} \left(\mathbf{E}\left[\left| \frac{X_i}{\sigma_n^*}\right|^p|\mathcal{G}_n^* \right] \right)^{1/p}, \text{ and  } \tilde \mu^*_p \equiv \max_{i \in N_n \backslash N_n^*} \left(\mathbf{E}\left[\left| \frac{X_i}{\sigma_n^*}\right|^p|\mathcal{G}_n^* \right] \right)^{1/p},
\end{eqnarray*}
where $\sigma_n^{*2} \equiv Var(\sum_{i \in N_n^*} X_i|\mathcal{G}_n^*)$, and 
\begin{eqnarray}
\label{Gn star}
	\mathcal{G}_n^* \equiv \mathcal{G} \vee \mathcal{M}_{N_n \backslash N_n^*}.
\end{eqnarray}
Note that the domain of the maximum in the definition of $\mu^*_p$ is $N_n^*$ whereas that of $\tilde \mu^*_p$ is $N_n\backslash N_n^*$. We also define for $i,j \in N_n^*$,
\begin{eqnarray*}
	r_n^{*2} \equiv \frac{1}{\sigma_n^{*4}}\mathbf{E}\left[ \left( \sum_{i \in N_n^*} \sum_{j \in \bar \nu_n(i) \cap N_n^*} \xi_{ij}\right)^2 |\mathcal{G}_n^*\right].
\end{eqnarray*}
\medskip

\begin{theorem}
	\label{thm: CLT-CHD}
	Suppose that a triangular array $\{X_i\}_{i \in N_n}$ is CND with respect to $(\nu_n, \mathcal{M})$ and $\mathbf{E}[X_i|\mathcal{M}_{\nu_n(i)}]=0$, a.e. for each $i \in N_n$ as in Theorem \ref{thm: CLT}, and let for $1\le r \le 4$ and $\varepsilon_n \ge 0$,
	\begin{eqnarray}
	\label{event 2}
	\mathcal{A}_{n,r}(\varepsilon_n) \equiv\{\omega \in \Omega: n^* d_{mx}^* d_{av}^* \mu_3^{*3} \le 1 \text{ and } (n- n^*) \tilde \mu_r^* + \rho_r^* \le \varepsilon_n\},
	\end{eqnarray}
	where
	\begin{eqnarray*}
		\rho_r^* \equiv \left( \mathbf{E}\left[ \left| \frac{1}{\sigma_n^*}\sum_{i \in N_n^*} \mathbf{E}[X_i|\mathcal{M}_{\nu_n^*(i)}^*] \right|^r \right] \right)^{1/r}.
	\end{eqnarray*}

	Then there exists an absolute constant $C>0$ such that on the event $\mathcal{A}_n(\varepsilon_n)$ with any constant $\varepsilon_n>0$, for each $t \in \mathbf{R}$, $n \ge 1$, and for each $1 \le r \le 4$, 
	\begin{eqnarray*}
		\sup_{t \in \mathbf{R}} \Delta^*(t;\mathcal{G}_n^*) &\le& C \sqrt{n^* d_{mx}^* d_{av}^* \mu_3^{*3}} - C \log(n^* d_{mx}^* d_{av}^*\mu_3^{*3}) \sqrt{n^* d_{mx}^{*2} d_{av}^*\mu_4^{*4} + r_n^{*2}} \\
		&+& C \varepsilon_n^{r/(r+1)}, \text{ a.e.},
	\end{eqnarray*} 
	where
	\begin{eqnarray*}
		\Delta^*(t;\mathcal{G}_n^*) \equiv \left|P\left\{\frac{1}{\sigma_n^*} \sum_{i \in N_n}  X_i \le t|\mathcal{G}_n^* \right\}  - \Phi(t)\right|.
	\end{eqnarray*}
\end{theorem}

Compared to Theorem \ref{thm: CLT}, the bound involves an additional term. This additional term arises because the sum may not be centered around zero when we condition on $\mathcal{M}_{\nu_n^*(i)}^*$. If $|X_i| \le 1/\sqrt{n}$ and $(\sigma_n^*)^2 \ge c$ for some $c>0$, we have $(n-n^*)\tilde \mu_r^* = O((n-n^*)/\sqrt{n})$. Furthermore, if $X$ has $\nu_n$ as a conditional dependency graph, we have $\rho_r^* = 0$ (because $\mathcal{M}_i = \mathcal{G}$ for all $i \in N_n$ in this case) and hence as long as
\begin{eqnarray*}
	(n - n^*)/\sqrt{n} \rightarrow 0,
\end{eqnarray*}
as $n \rightarrow \infty$, the third term in the bound vanishes, and the same CLT as in Theorem \ref{CLT cond} is restored. In this case, if $(n^*/n)^{1/2}(((n-n^*)/\sqrt{n})^{4/5} + (d_{mx}^*d_{av}^*/\sqrt{n^*})^{1/2})$ converges to zero faster than $(d_{mx}d_{av}/\sqrt{n})^{1/2}$, Theorem \ref{thm: CLT-CHD} has an improved rate over Theorem \ref{thm: CLT}. Such an approximation captures the situation where the neighborhood system $\nu_n$ has a very small fraction of very high degree vertices.

Such an improvement can still arise generally, even if $X$ does not have $\nu_n$ as a conditional dependency graph. To see this, note that for each $i \in N_n \backslash \bar \nu_n(N_n \backslash N_n^*)$, $X_i$ is conditionally independent of $\mathcal{M}_{N_n \backslash N_n^*}$ given $\mathcal{M}_{\nu_n(i)}$, we have
\begin{eqnarray*}
	\mathbf{E}[X_i|\mathcal{M}_{\nu_n^*(i)}^*] &=& \mathbf{E}[X_i|\mathcal{M}_{(\nu_n(i) \cap N_n^*)\cup (N_n \backslash N_n^*)}]\\
	&=& \mathbf{E}[X_i|\mathcal{M}_{\nu_n(i) \cup (N_n \backslash N_n^*)}] = \mathbf{E}[X_i|\mathcal{M}_{\nu_n(i)}] = 0.
\end{eqnarray*}
Hence 
\begin{eqnarray*}
	\rho_r^* = \left( \mathbf{E}\left[ \left| \frac{1}{\sigma_n^*}\sum_{i \in \bar \nu_n(N_n \backslash N_n^*) \cap N_n^*} \mathbf{E}[X_i|\mathcal{M}_{\nu_n^*(i)}^*] \right|^r \right] \right)^{1/r}.
\end{eqnarray*}
Suppose that $|X_i| \le 1/\sqrt{n}$, $(\sigma_n^*)^2 \ge c$ for some $c>0$, $n - n^* = O(1)$, and $d_{av} = O(1)$ and $d_{av}^* = O(1)$. Then we have
\begin{eqnarray}
\label{rho star}
	\rho_r^* = O(|\bar \nu_n(N_n \backslash N_n^*)|/\sqrt{n}) = O(d_{mx}/\sqrt{n}).
\end{eqnarray}
The rate in Theorem \ref{thm: CLT-CHD} improves on that in Theorem \ref{thm: CLT} because with $p>1$, $(\rho_p^*)^{4/5} = O((d_{mx}/\sqrt{n})^{4/5}) = o((d_{mx}/\sqrt{n})^{1/2})$.

As we shall see later in Section \ref{subsubsec: CND cond on high deg}, the approach of CND conditional on high-degree vertices is useful for obtaining stable central limit theorem for empirical processes when the random variables are CND with respect to a neighborhood system having a maximum degree $d_{mx}$ increasing to infinity as $n \rightarrow \infty$.

\subsection{Empirical Processes and Stable Convergence}
\subsubsection{Stable Convergence in Metric Spaces}
Let us first introduce some preliminary results about stable convergence in metric spaces. Let $(\mathbb{D},d)$ be a given metric space and define $\mathcal{B}(\mathbb{D})$ to be the Borel $\sigma$-field of $\mathbb{D}$. Let $(\Omega,\mathcal{F},P)$ be the given probability space, where $\mathcal{G}$ is a sub $\sigma$-field of $\mathcal{F}$. 
Recall that a map $L:\mathcal{B}(\mathbb{D}) \times \Omega \rightarrow [0,1]$ called a \textit{Markov kernel} if for each $\omega \in \Omega$, $L(\cdot,\omega)$ is a Borel probability measure on $\mathcal{B}(\mathbb{D})$ and for each $B \in \mathcal{B}(\mathbb{D})$, $L(B,\cdot)$ is $\mathcal{F}$-measurable. Following the notation in \cite{Hausler/Luschgy:10:PTSM}, let us define the \textit{marginal} of $L$ on $\mathcal{B}(\mathbb{D})$ as
\begin{eqnarray*}
	PL(B) \equiv \int L(B,\omega)dP(\omega).
\end{eqnarray*}

For a given finite collection $\{h_1,...,h_k\} \subset \mathcal{H}$, the finite dimensional projection of a Markov kernel $L$ is defined to be a Markov kernel $L_{h_1,...,h_k}:\mathcal{B}(\mathbf{R}^k) \times \Omega \rightarrow [0,1]$ such that for any $B \in \mathcal{B}(\mathbf{R}^k)$,
\begin{eqnarray*}
	L_{h_1,...,h_k}(B,\omega) \equiv L(\{y \in \mathbb{D}:(y(h_1),...,y(h_k)) \in  B\},\omega).	
\end{eqnarray*}
In the spirit of the Hoffman-Jorgensen approach, we consider the following definition of stable convergence for empirical processes. (See \cite{Berti/Pratelli/Rigo:12:EJP}, p.2., for a similar definition.)

\begin{definition}
	\label{def: Stable Conv}
	\noindent Suppose that we are given a sub $\sigma$-field $\mathcal{G}$, a sequence of $\mathbb{D}$-valued stochastic processes $\zeta_n$, a Markov kernel $L$ on $\mathcal{B}(\mathbb{D}) \times \Omega$, and a $\mathbb{D}$-valued Borel measurable random element $\zeta$ that has a Markov kernel $L$.
	
	Suppose that for each $U \in \mathcal{G}$, and each bounded Lipschitz functional $f$ on $\mathbb{D}$,
	\begin{eqnarray*}
		\mathbf{E}^*[f(\zeta_n)1_U] \textnormal{  } \rightarrow  \textnormal{  } \int_U \int_\mathbb{D} f(y) L(dy,\omega)dP(\omega),
	\end{eqnarray*}
	as $n \rightarrow \infty$, where $\mathbf{E}^*$ denotes the outer-expectation and $1_U$ is the indicator function of the event $U$. Then we say that $\zeta_n$ \textit{converges to} $L$, $\mathcal{G}$\textit{-stably}, (or equivalently, $\zeta_n$ \textit{converges to} $\zeta$, $\mathcal{G}$-\textit{stably}), and write
	\begin{eqnarray*}
		\zeta_n \rightarrow L, \mathcal{G}\text{-stably } \text{(or equivalently, } \zeta_n \rightarrow \zeta, \mathcal{G}\text{-stably.)}
	\end{eqnarray*}
\end{definition}

Stable convergence according to Definition \ref{def: Stable Conv} implies weak convergence (in the sense of Hoffman-Jorgensen). When $\zeta_n$ is Borel measurable, the above definition is equivalent to the weak convergence of Markov kernels and many of the existing results on stable convergence carry over. However, this equivalence does not extend to the case of $\zeta_n$ being non-measurable, because there is no proper definition of Markov kernels for nonmeasurable stochastic processes. Nevertheless the above definition can still be useful when one needs to deal with random norming, as shown in the following lemma which generalizes Theorem 1' in \cite{Aldous/Eagleson:78:AP} and part of Theorem 3.18 of \cite{Hausler/Luschgy:10:PTSM}.

\begin{lemma}
	\label{generalized slutsky}
	Suppose that $\zeta_n,\zeta,\xi_n,\xi$ are $\mathbb{D}$-valued random variables, and that $\zeta_n \rightarrow \zeta$, $\mathcal{G}$-stably, where $\zeta$ and $\xi$ are Borel measurable, and let $P_\zeta$ denote the distribution of $\zeta$. Then the following holds.
	\medskip
	
	\noindent (i) If $P^*\{d(\xi_n,\xi) > \varepsilon\} \rightarrow 0$ as $n \rightarrow \infty$ for each $\varepsilon>0$, and $\xi$ is $\mathcal{G}$-measurable, then
	\begin{eqnarray*}
		(\zeta_n,\xi_n) \rightarrow (\zeta,\xi), \mathcal{G}\text{-stably},
	\end{eqnarray*}
    where $P^*$ denotes the outer probability.
	
	\noindent (ii) If $f:\mathbb{D} \rightarrow \mathbb{D}$ is $P_\zeta$-a.e. continuous, then $f(\zeta_n) \rightarrow f(\zeta)$, $\mathcal{G}$-stably.
\end{lemma}

The first result is a stable-convergence analogue of Cram\'{e}r-Slutsky lemma. The second result is a continuous mapping theorem.

\subsubsection{Stable Convergence of an Empirical Process}
Suppose that $\{Y_{i}\}_{i \in N_n}$ is a given triangular array of $\mathbf{R}$-valued random variables which is CND with respect to $(\nu_n,\mathcal{M})$.  Let $\mathcal{H}$ be a given class of real measurable functions on $\mathbf{R}$, having a measurable envelope $H$. Then, we consider the following empirical process:
\begin{eqnarray*}
	\{ \mathbb{G}_n(h):h \in \mathcal{H} \},
\end{eqnarray*}
where, for each $h \in \mathcal{H}$,
\begin{eqnarray*}
	\label{nu g}
	\mathbb{G}_n(h)\equiv \frac{1}{\sqrt{n}} \sum_{i \in N_n} (h(Y_{i}) - \mathbf{E}[h(Y_{i})|\mathcal{M}_{\nu_n(i)}]).
\end{eqnarray*}
The empirical process $\nu_n$ takes a value in $l^\infty(\mathcal{H})$, the collection of bounded functions on $\mathcal{H}$ which is endowed with the sup norm so that $(l^\infty(\mathcal{H}),\|\cdot\|_\infty)$ forms the metric space $(\mathbb{D},d)$ with $
	\| h \|_\infty \equiv \sup_{y \in \mathbf{R}} |h(y)|.$
In this section, we explore conditions for the class $\mathcal{H}$ and the joint distribution of the triangular array $\{Y_{i}\}_{i \in N_n}$ which delivers the stable convergence of the empirical process. Stable convergence in complete separable metric spaces can be defined as a weak convergence of Markov kernels. (See \cite{Hausler/Luschgy:10:PTSM}.) However, this definition does not extend to the case of empirical processes taking values in $\mathbb{D}$ that is endowed with the sup norm, due to non-measurability. 

Weak convergence of an empirical process to a Gaussian process is often established in three steps. First, we show that the class of functions is totally bounded with respect to a certain pseudo-metric $\rho$. Second, we show that each finite dimensional projection of the empirical process converges in distribution to a multivariate normal random vector. Third, we establish the asymptotic $\rho$-equicontinuity of the empirical process.

Let $\rho$ be a given pseudo-metric on $\mathcal{H}$ such that $(\mathcal{H},\rho)$ is a totally bounded metric space. Then we define
\begin{eqnarray*}
	U_\rho (\mathcal{H}) \equiv \{ y \in \mathbb{D}: y \text{ is uniformly } \rho\text{-continuous on } \mathcal{H}\}.
\end{eqnarray*}
The following theorem shows that we can take a similar strategy in proving the stable convergence of an empirical process to a Markov kernel. The structure and the proof of the theorem is adapted from Theorem 10.2 of \cite{Pollard:90:EmpiricalProcesses}.
\begin{theorem}
	\label{thm: fundamental theorem}
	Suppose that the stochastic process $\zeta_n \in \mathbb{D}$ is given, where $(\mathcal{H},\rho)$ is a totally bounded metric space. Suppose that the following conditions hold.
	
	\noindent (i) For each finite set $\{h_1,...,h_k\} \subset \mathcal{H}$, $(\zeta_n(h_1),...,\zeta_n(h_k)) \rightarrow L_{h_1,...,h_k}$, $\mathcal{G}$-stably, where $L_{h_1,...,h_k}$ is a Markov kernel on $\mathcal{B}(\mathbf{R}^k) \times \Omega$.
	
	\noindent (ii) $\zeta_n$ is asymptotically $\rho$-equicontinuous on $\mathcal{H}$, i.e., for each $\varepsilon>0,\eta>0$, there exists $\delta>0$ such that
	\begin{eqnarray*}
		\textnormal{limsup}_{n \rightarrow \infty} P^*\left\{\sup_{h,h' \in \mathcal{H}:\rho(h,h')<\delta} |\zeta_n(h) - \zeta_n(h')| > \eta \right\} < \varepsilon,
	\end{eqnarray*}
    where $P^*$ denotes the outer probability.
	
	Then there exists a Markov kernel $L$ on $\mathcal{B}(\mathbb{D}) \times \Omega$ such that the following properties are satisfied.
	
	(a) The finite dimensional projections of $L$ are given by Markov kernels $L_{h_1,...,h_k}$.
	
	(b) $PL(U_\rho(\mathcal{H}))=1$.
	
	(c) $\zeta_n \rightarrow L$, $\mathcal{G}$-stably.
	
	Conversely, if $\zeta_n$ converges to Markov kernel $L$ on $\mathcal{B}(\mathbb{D}) \times \Omega$, $\mathcal{G}$-stably, where $PL(U_\rho(\mathcal{H}))=1$, then (i) and (ii) are satisfied.
\end{theorem}

It is worth noting that for stable convergence of empirical processes, the conditions for the asymptotic equicontinuity and the totally boundedness of the function class with respect to a pseudo-norm are as in the standard literature on weak convergence of empirical processes. The only difference here is that the convergence of finite dimensional distributions is now replaced by the stable convergence of finite dimensional projections.

\subsubsection{Maximal Inequality}
This subsection presents a maximal inequality in terms of bracketing entropy bounds. The maximal inequality is useful primarily for establishing asymptotic $\rho$-equicontinuity of the empirical process but also for many other purposes. First, we begin with a tail bound for a sum of CND random variables. \cite{Janson:04:RSA} established an exponential bound for a sum of random variables that have a dependency graph. The following exponential tail bound is crucial for our maximal inequality. The result below is obtained by slightly modifying the proof of Theorem 2.3 of \cite{Janson:04:RSA}.

\begin{lemma}
	\label{lemma: tail bound}
	\noindent Suppose that $\{X_{i}\}_{i \in N_n}$ is a triangular array of random variables that take values in $[-M,M]$ and are CND with respect to $(\nu_n,\mathcal{M})$, with $\mathbf{E}[X_i|\mathcal{M}_{\nu_n(i)}]=0$, and let $\sigma_i^2 \equiv Var\left(X_i|\mathcal{M}_{\nu_n(i)}\right)$ and $V_n \equiv \sum_{i \in N_n} \mathbf{E}[\sigma_i^2|\mathcal{G}]$ with $\mathcal{G}$ as defined in (\ref{sigma G}).
	
	Then, for any $\eta > 0$, 
	\begin{eqnarray}
	\label{bd2} \quad \quad
	P\left\{ \left|\sum_{i \in N_n} X_i \right| \ge \eta |\mathcal{G} \right\}
	&\le& 2 \exp \left(-\frac{\eta^2}{2 (d_{mx} + 1) (2(d_{mx} + 1) V_n + M\eta/3)} \right), \text{a.e.},
	\end{eqnarray}
	for all $n \ge 1$.
	
	Furthermore, if Condition A holds and the condition $\mathbf{E}[X_{i}|\mathcal{M}_{\nu_n(i)}] = 0$ is replaced by $\mathbf{E}[X_{i}|\mathcal{G}] = 0$ and the $\sigma$-fields $\mathcal{M}_{\nu_n(i)}$ in $\sigma_i$'s are replaced by $\mathcal{G}$, then the following holds: for any $\eta > 0$, 
	\begin{eqnarray}
	\label{bd3} \quad \quad
	P\left\{ \left|\sum_{i \in N_n} X_i \right| \ge \eta |\mathcal{G} \right\}
	\le 2 \exp \left(-\frac{8 \eta^2}{25 (d_{mx} + 1) ( V_n +  M \eta/3)} \right), \text{a.e.},
	\end{eqnarray}
	for all $n \ge 1$.
\end{lemma}

The bound in (\ref{bd3}) is the one obtained by \cite{Janson:04:RSA} for the case of dependency graphs. From this, the following form of maximal inequality for a finite set immediately follows from Lemma A.1 of \cite{vanderVaart:96:AS}.

\begin{corollary}
	\label{maximal inequality0}
	\noindent Suppose that $\{Y_i\}_{i \in N_n}$ is a triangular array of random variables that are CND with respect to $(\nu_n,\mathcal{M})$. Let for each $h \in \mathcal{H}$, $V_n(h) \equiv n^{-1}\sum_{i \in N_n} \mathbf{E}[\sigma_i^2(h)|\mathcal{G}]$ and $\sigma_i^2(h) \equiv Var(h(Y_i)|\mathcal{M}_{\nu_n(i)})$ with $\mathcal{G}$ as defined in (\ref{sigma G}).
	
	Then there exists an absolute constant $C>0$ such that
	\begin{eqnarray*}
		&& \mathbf{E}\left[ \max_{1 \le s \le m} |\mathbb{G}_n(h_s)||\mathcal{G}\right] \\ &\le& C (d_{mx} +1) \left(\frac{ J}{\sqrt{n}} \log(1+m) + \sqrt{\log(1+m)\max_{1\le s \le m} V_n(h_s)}\right), \text{a.e.},
	\end{eqnarray*}
	for any $n\ge 1$ and any $m \ge 1$ with a finite subset $\{h_1,...,h_m\}$ of $\mathcal{H}$ such that for some constant $J>0$, $\max_{1 \le s \le m}\sup_{x \in \mathbf{R}}|h_s(x)| \le J$.
\end{corollary}
\medskip

Let us now elevate the above inequality to a maximal inequality over function class $\mathcal{H}$. Recall that we allow the random variables $Y_i$'s to be idiosyncratically distributed across $i$'s. 

We define the following semi-norm on $\mathcal{H}$:
\begin{eqnarray*}
	\bar \rho_n (h)  &\equiv&   \sqrt{\frac{1}{n} \sum_{i \in N_n} \mathbf{E}[h^2(Y_i)] }.
\end{eqnarray*} 
We denote $N_{[]}(\varepsilon,\mathcal{H},\bar \rho_n)$ to be the $\varepsilon$-bracketing number of $\mathcal{H}$ with respect to $\bar \rho_n$, i.e., the smallest number $J$ of the brackets $[h_{L,j}, h_{U,j}]$, $j=1,...,J$, such that $\bar \rho_n(h_{U,j} - h_{L,j}) \le \varepsilon$. The following lemma establishes the maximal inequality in terms of a bracketing entropy bound.

\begin{lemma} [Maximal Inequality]
	\label{lemm: maximal inequality}
    \noindent Suppose that $\{Y_{i}\}_{i\in N_n}$ is a
	triangular array of random variables that are CND with respect to $(\nu_n,
	\mathcal{M})$. Suppose further that the class $\mathcal{H}$ of functions
	have an envelope $H$ such that $\bar{\rho}_{n}(H)<\infty $. Then, there
	exists an absolute constant $C>0$ such that for each $n \ge 1$,
	\begin{equation*}
	\mathbf{E}^*\left[ \sup_{h \in \mathcal{H}} |\mathbb{G}_n(h)|\right] \leq
	C(1+d_{mx})\int_{0}^{\bar{\rho}_{n}(H)}\sqrt{1+\log N_{[]}(\varepsilon ,\mathcal{H},\bar{\rho}_{n})}d\varepsilon.
	\end{equation*}
\end{lemma}
\medskip

The bracketing entropy bound in Lemma \ref{lemm: maximal inequality} involves the maximum degree $d_{mx}$. Hence the bound is useful only when the neighborhood system $\nu_n$ does not have a maximum degree increasing with $n$.

\subsubsection{Stable Central Limit Theorem}

First, let us say that a stochastic process $\{\mathbb{G}(h): h \in \mathcal{H}\}$ is a \textit{$\mathcal{G}$-mixture Gaussian process} if for any finite collection $\{h_1,...,h_m\} \subset \mathcal{H}$, the distribution of random vector $[\mathbb{G}(h_1),\mathbb{G}(h_2),...,\mathbb{G}(h_m)]$ conditional on $\mathcal{G}$ is a multivariate normal distribution. Also, we call a Markov kernel $\mathcal{K}$ a \textit{$\mathcal{G}$-mixture Gaussian Markov kernel} associated with a given $\mathcal{G}$-mixture Gaussian process $\mathbb{G}$ if for any $h_1,...,h_m \in \mathcal{H}$, the conditional distribution of $[\mathbb{G}(h_1),\mathbb{G}(h_2),...,\mathbb{G}(h_m)]$ given $\mathcal{G}$ is given by the finite dimensional projection $\mathcal{K}_{h_1,...,h_m}$ of $\mathcal{K}$. Let us summarize the conditions as follows.

\begin{assumption}
	\label{assump: conds emp CLT}
\noindent (a) There exists $C>0$ such that for all $n \ge 1$,
	\begin{eqnarray*}
		(1+d_{mx})\int_{0}^{\bar{\rho}_{n}(H)}\sqrt{1+\log N_{[]}(\varepsilon ,%
			\mathcal{H},\bar{\rho}_{n})}d\varepsilon < C.
	\end{eqnarray*}

\noindent (b) For any $h_1,h_2 \in \mathcal{H}$,
\begin{eqnarray*}
	\mathbf{E}[\mathbb{G}_n(h_1)\mathbb{G}_n(h_2)|\mathcal{G}] \rightarrow_P K(h_1,h_2|\mathcal{G}), 
\end{eqnarray*}
for some $K(\cdot,\cdot|\mathcal{G})(\omega):\mathcal{H}\times\mathcal{H} \rightarrow \mathbf{R}, \omega \in \Omega$, which is positive semidefinite a.e., and non-constant at zero. 
\medskip

\noindent (c) For each $h \in \mathcal{H}$, $- n^{-1/2} r_n(h)\log(n^{-1/2} d_{mx} d_{av})  \rightarrow_P 0$, as $n \rightarrow \infty$,
	where
	\begin{eqnarray*}
		r_n^2(h) &\equiv& \mathbf{E}\left[\left(\frac{1}{\sqrt{n}}\sum_{i \in N_n}\sum_{j \in \bar \nu_n(i)} \xi_{ij}(h)\right)^2|\mathcal{G}\right],
	\end{eqnarray*}
	with
	\begin{eqnarray*}
		\xi_{ij}(h) &\equiv& \mathbf{E}[(h(X_i)-\mathbf{E}[h(X_i)|\mathcal{M}_{\nu_n(i)}])
		(h(X_j)-\mathbf{E}[h(X_j)|\mathcal{M}_{\nu_n(j)}])|\mathcal{M}_{\nu_n(i,j)}]\\
		&& - \mathbf{E}[(h(X_i)-\mathbf{E}[h(X_i)|\mathcal{M}_{\nu_n(i)}])
		(h(X_j)-\mathbf{E}[h(X_j)|\mathcal{M}_{\nu_n(j)}])|\mathcal{G}].
	\end{eqnarray*}
\medskip

\noindent (d) For each $h \in \mathcal{H}$, $\rho(h) \equiv \lim_{n \rightarrow \infty} \bar \rho_n(h)$ exists in $[0,\infty)$, and satisfies that whenever $\rho(h_n) \rightarrow 0$ as $n \rightarrow \infty$, $\bar \rho_n(h_n) \rightarrow 0$ as $n \rightarrow \infty$ as well.
\end{assumption}

The following result gives a Donsker-type stable convergence of empirical processes.

\begin{theorem}
	\label{thm: empirical CLT}
	Suppose that $\{Y_{i}\}_{i\in N_n}$ is a
	triangular array of random variables that are CND with respect to $(\nu_n,
	\mathcal{M})$, satisfying Assumption \ref{assump: conds emp CLT}. Suppose further that there exists $C>0$ such that for each $n \ge 1$, $\max_{i \in N_n} \mathbf{E}[H(Y_{i})^4] < C$, where $H$ is an envelope of $\mathcal{H}$.
	
	Then $\nu_n$ converges to a $\mathcal{G}$-mixture Gaussian process $\mathbb{G}$ in $l^\infty(\mathcal{H})$, $\mathcal{G}$-stably, such that for any $h_1,h_2 \in \mathcal{H}$,
		$\mathbf{E}[\mathbb{G}(h_1)\mathbb{G}(h_2)|\mathcal{G}] = K(h_1,h_2|\mathcal{G}), a.e.$
		
	Furthermore, we have $P\mathcal{K}(U_\rho(\mathcal{H}))=1$, where $\mathcal{K}$ is the $\mathcal{G}$-mixture Gaussian Markov kernel associated with $\mathbb{G}$.
\end{theorem}

The fourth moment condition $\max_{i \in N_n} \mathbf{E}[H(Y_{i})^4] < \infty$ is used to ensure the convergence of finite dimensional distributions using Theorem \ref{thm: CLT}. It is worth noting that Assumption \ref{assump: conds emp CLT}(a) essentially requires that the maximum degree $d_{mx}$ to be bounded. It is interesting that this condition was not required for the CLT in Theorem \ref{thm: CLT}. This stronger condition for the maximum degree is used to establish the asymptotic equicontinuity of the process $\{\mathbb{G}_n(h): h \in \mathcal{H}\}$. When the neighborhood system $\nu_n$ is generated according to a model of stochastic graph formation, this condition is violated for many existing models of graph formation used, for example, for social network modeling. In the next section, we utilize the approach of conditioning on high-degree vertices to weaken this condition.
 
\subsubsection{Conditional Neighborhood Dependence Conditional on High-Degree Vertices}
\label{subsubsec: CND cond on high deg}
As mentioned before, Assumption \ref{assump: conds emp CLT}(a) requires that $d_{mx}$ be bounded. Following the idea of conditioning on high degree vertices as in Theorem \ref{thm: CLT-CHD}, let us explore a stable convergence theorem that relaxes this requirement. As we did prior to Theorem \ref{thm: CLT-CHD}, we choose $N_n^* \subset N_n$ to be a given subset and let $d_{mx}^*$ and $d_{av}^*$ be as defined in (\ref{max deg star}).

First, write
\begin{eqnarray}
\label{decomp}
	\mathbb{G}_n(h) = \mathbb{G}_n^*(h)  + R_n^*(h) + \rho_n^*(h),
\end{eqnarray}
where
\begin{eqnarray*}
	\mathbb{G}_n^*(h) &\equiv& \frac{1}{\sqrt{n}} \sum_{i \in N_n^*} \left(h(Y_i) - \mathbf{E}[h(Y_i)|\mathcal{M}_{\nu_n^*(i)}^*] \right) \\
	R_n^*(h) &\equiv& \frac{1}{\sqrt{n}} \sum_{i \in N_n \backslash N_n^*} \left(h(Y_i) - \mathbf{E}[h(Y_i)|\mathcal{M}_{\nu_n(i)}] \right),
	 \text{ and } \\
	\rho_n^*(h) &\equiv& \frac{1}{\sqrt{n}} \sum_{i \in N_n^*} \left(\mathbf{E}[h(Y_i)|\mathcal{M}_{\nu_n(i)}] - \mathbf{E}[h(Y_i)|\mathcal{M}^*_{\nu_n^*(i)}] \right). 
\end{eqnarray*}
Note that
\begin{eqnarray*}
	\mathbf{E}^*\left[\sup_{h \in \mathcal{H}}|R_n^*(h)|\right]
	&\le& \frac{2(n - n^*)}{\sqrt{n}} \max_{i \in N_n} \sqrt{\mathbf{E}[H^2(Y_i)]}.
\end{eqnarray*}
Since $\{Y_i\}_{i \in N_n^*}$ is CND with respect to $(\nu_n^*,\mathcal{M}^*)$ as defined in (\ref{nu star}) and (\ref{M star}), we can apply the previous results to $\mathbb{G}_n^*(h)$. This gives the following extension of the maximal inequality in Lemma \ref{lemm: maximal inequality}. Since the maximal inequality is often of independent interest, let us state it formally.

\begin{lemma} [Maximal Inequality]
	\label{maximal inequality 2}
	\noindent Suppose that $\{Y_{i}\}_{i\in N_n}$ is a
	triangular array of random variables that are CND with respect to $(\nu_n,
	\mathcal{M})$. Suppose further that the class $\mathcal{H}$ of functions
	have an envelope $H$ such that $\bar{\rho}_{n}(H)<\infty $.
	
	Then, there exists an absolute constant $C>0$ such that for all $n \ge 1$,
	\begin{eqnarray*}
	\mathbf{E}^*\left[ \sup_{h \in \mathcal{H}} |\mathbb{G}_n(h)|\right] &\le&
	\frac{C\sqrt{n^*}(1+d_{mx}^*)}{\sqrt{n}}\int_{0}^{\bar{\rho}_{n}(H)}\sqrt{1+\log N_{[]}(\varepsilon,
		\mathcal{H},\bar{\rho}_{n})}d\varepsilon \\
	&+& C \left(\frac{n - n^*}{\sqrt{n}} \max_{i \in N_n} \sqrt{\mathbf{E}[H^2(Y_i)]} + \mathbf{E}^*\left[\sup_{h \in \mathcal{H}} \rho_n^*(h) \right]\right).
	\end{eqnarray*}
\end{lemma}
\medskip
If $\| H\|_\infty < C$ and $n - n^* = O(1)$, the second term in the bound is $O(d_{mx}/\sqrt{n})$ similarly as we derived in (\ref{rho star}). Thus this bound is an improvement over Lemma \ref{lemm: maximal inequality}, whenever $O(d_{mx}^*) =  o (d_{mx})$ as $n \rightarrow \infty$. Let us turn to the Donsker-type stable convergence of an empirical process. We modify Assumption \ref{assump: conds emp CLT} as follows.
\begin{assumption}
	\label{assump: conds emp CLT 2}
\noindent (a) There exists $C>0$ such that for all $n \ge 1$,
	\begin{eqnarray*}
		(1+d_{mx}^*)\int_{0}^{\bar{\rho}_{n}(H)}\sqrt{1+\log N_{[]}(\varepsilon ,%
			\mathcal{H},\bar{\rho}_{n})}d\varepsilon < C.
	\end{eqnarray*}

\noindent (b) For any $h_1,h_2 \in \mathcal{H}$,
\begin{eqnarray*}
	\mathbf{E}[\mathbb{G}_n(h_1)\mathbb{G}_n(h_2)|\mathcal{G}_n^*] \rightarrow_P K(h_1,h_2|\mathcal{G}^*), 
\end{eqnarray*}
for some $K(\cdot,\cdot|\mathcal{G}^*)(\omega):\mathcal{H}\times\mathcal{H} \rightarrow \mathbf{R}, \omega \in \Omega$, which is positive semidefinite a.e, and non-constant at zero, and for some sub $\sigma$-field $\mathcal{G}^*$ of $\mathcal{F}$, where $\mathcal{G}_n^*$ is as defined in (\ref{Gn star}).
\medskip

\noindent (c) For each $h \in \mathcal{H}$, $- n^{*-1/2} r_n^*(h)\log(n^{*-1/2} d_{mx}^* d_{av}^*)  \rightarrow_P 0$, as $n \rightarrow \infty$,
where
\begin{eqnarray*}
	r_n^{*2}(h) &\equiv& \mathbf{E}\left[\left(\frac{1}{\sqrt{n}}\sum_{i \in N_n^*}\sum_{j \in \bar \nu_n(i) \cap N_n^*} \xi_{ij}(h)\right)^2|\mathcal{G}_n^*\right].
\end{eqnarray*}

\noindent (d) For each $h \in \mathcal{H}$, $\rho(h) \equiv \lim_{n \rightarrow \infty} \bar \rho_n(h)$ exists in $[0,\infty)$ and satisfies that whenever $\rho(h_n) \rightarrow 0$ as $n \rightarrow \infty$, $\bar \rho_n(h_n) \rightarrow 0$ as $n \rightarrow \infty$ as well.\medskip

\noindent (e) $\mathbf{E}^*\left[\sup_{h \in \mathcal{H}} \rho_n^*(h) \right] \rightarrow 0$, as $n \rightarrow \infty$. \medskip
\end{assumption}

While Condition (a) essentially requires that $d_{mx}^*$ be bounded, Condition (c) allows $d_{mx}$ to increase to infinity as $n \rightarrow \infty$. The condition in (b) that $K(\cdot,\cdot|\mathcal{G}^*)$ be non-constant at zero requires that
\begin{eqnarray*}
	n - n^* \rightarrow \lambda, \text{ as } n \rightarrow \infty,
\end{eqnarray*}
for some $\lambda \in [0, \infty)$. Thus the number of the high degree vertices $(n-n^*)$ selected when we set $N_n^* \subset N_n$ should be bounded as $n \rightarrow \infty$. In combination with (e), this implies that we have $(n-n^*) / \sqrt{n} + \mathbf{E}^*\left[\sup_{h \in \mathcal{H}} \rho_n^*(h) \right] \rightarrow 0$ as $n \rightarrow \infty$, which makes it suffice to focus on $\mathbb{G}_n^*(h)$ in (\ref{decomp}) for a stable limit theorem. We obtain the following extended version of Theorem \ref{thm: empirical CLT}.

\begin{theorem}
	\label{thm: empirical CLT 2}
	Suppose that $\{Y_{i}\}_{i\in N_n}$ is a
	triangular array of random variables that are CND with respect to $(\nu_n,
	\mathcal{M})$, satisfying Assumption \ref{assump: conds emp CLT 2}. Suppose further that there exists $C>0$ such that for each $n \ge 1$, $\max_{i \in N_n} \mathbf{E}[H(Y_{i})^4] < C$, where $H$ is an envelope of $\mathcal{H}$.
	
	Then $\nu_n$ converges to a $\mathcal{G}^*$-mixture Gaussian process $\mathbb{G}^*$ in $l^\infty(\mathcal{H})$, $\mathcal{G}^*$-stably, such that for any $h_1,h_2 \in \mathcal{H}$, $
		\mathbf{E}[\mathbb{G}^*(h_1)\mathbb{G}^*(h_2)|\mathcal{G}^*] = K(h_1,h_2|\mathcal{G}^*), a.e.$
	
	Furthermore, we have $P\mathcal{K}(U_\rho(\mathcal{H}))=1$, where $\mathcal{K}$ is the $\mathcal{G}^*$-mixture Gaussian Markov kernel associated with $\mathbb{G}^*$.
\end{theorem}

If we take $N_n^*$ to be identical to $N_n$, Theorem \ref{thm: empirical CLT 2} is reduced to Theorem \ref{thm: empirical CLT}. However, Theorem \ref{thm: empirical CLT 2} shows that approximation of the distribution of an empirical process by a mixture Gaussian process is possible even if $d_{mx} \rightarrow \infty$ as $n \rightarrow \infty$.

\section{Appendix: Mathematical Proofs}

To simplify the notation, we follow \cite{vanderVaart:96:AS} and write $a_n \lesssim b_n$ for any sequence of numbers, whenever $a_n \le C b_n$ for all $n \ge 1$ with some absolute constant $C>0$. The absolute constant can differ across different instances of $\lesssim$.

For any positive integer $k$ and $\mathbf{i} \equiv (i_1,...,i_k) \in N_n^k$ and a triangular array of random variables $\{X_i\}_{i \in N_n}$, we define
\begin{eqnarray}
    \label{X}
	X(\mathbf{i}) \equiv \prod_{r=1}^k X_{i_r}.
\end{eqnarray}
The following lemma is useful for the proofs of various results.  
\medskip

\begin{lemma}
	\label{lemm: eq4}
	Suppose that $\nu_n$ is a neighborhood system on $N_n$ and $\{X_i\}_{i \in N_n}$ is a triangular array of random variables that are CND with respect to $(\nu_n,\mathcal{M})$, where $\mathcal{M}=\{\mathcal{M}_i\}_{i \in N_n}$. Furthermore, for given positive integer $k$, let $\mathbf{i} = (i_1,...,i_k) \in N_n^k$ be such that it has two partitioning subvectors $I_{k,1}(\mathbf{i})$ and $I_{k,2}(\mathbf{i})$ of $\mathbf{i}$ such that the entries of $I_{k,1}(\mathbf{i})$ are from $N_n \setminus \bar \nu_n(I_{k,2}(\mathbf{i}))$ and the entries of $I_{k,2}(\mathbf{i})$ are from $N_n \setminus \bar \nu_n(I_{k,1}(\mathbf{i}))$. Then,
	\begin{eqnarray*}
		\mathbf{E}[X(\mathbf{i}) |\mathcal{M}_{ \nu_n(\mathbf{i})}]
		&=&\mathbf{E}[X(I_{k,1}(\mathbf{i})) |\mathcal{M}_{ \nu_n(I_{k,1}(\mathbf{i}))}]
		\mathbf{E}[X(I_{k,2}(\mathbf{i})) |\mathcal{M}_{ \nu_n(I_{k,2}(\mathbf{i}))}].
	\end{eqnarray*}
	
	Suppose further that Condition A holds. Then,
	\begin{eqnarray*}
		\mathbf{E}[X(\mathbf{i}) |\mathcal{G}]
		=\mathbf{E}[X(I_{k,1}(\mathbf{i})) |\mathcal{G}]
		\mathbf{E}[X(I_{k,2}(\mathbf{i})) |\mathcal{G}].
	\end{eqnarray*}
\end{lemma}

\noindent \textbf{Proof:} By the choice of $\mathbf{i}$, we have
\begin{eqnarray}
\label{nu}
	\nu_n(\mathbf{i}) &=& \nu_n(I_{k,1}(\mathbf{i})) \cup \nu_n(I_{k,2}(\mathbf{i})), \\ \notag
	\nu_n(I_{k,1}(\mathbf{i})) &\subset& \nu_n(\mathbf{i}), \text{ and } \nu_n(I_{k,2}(\mathbf{i})) \subset \nu_n(\mathbf{i}).
\end{eqnarray}
To see the second statement, note that whenever $i \in \nu_n(I_{k,1}(\mathbf{i}))$, we have $i \notin I_{k,1}(\mathbf{i})$ and $i \notin I_{k,2}(\mathbf{i})$, and there must exist $j \in I_{k,1}(\mathbf{i})$ such that $i \in \nu_n(j)$. Since $I_{k,1}(\mathbf{i}) \subset \mathbf{i}$, we find that $i \in \nu_n(\mathbf{i})$.

As for the first statement of the lemma, we write
\begin{eqnarray*}
	\label{eq22}
	\mathbf{E}[X(\mathbf{i}) |\mathcal{M}_{ \nu_n(\mathbf{i})}] 
	&=& \mathbf{E}[\mathbf{E}[X(I_{k,1}(\mathbf{i})) |\mathcal{M}_{ \nu_n(\mathbf{i})},X(I_{k,2}(\mathbf{i}))] X(I_{k,2}(\mathbf{i})) |\mathcal{M}_{ \nu_n(\mathbf{i})}]\\ \notag
	&=& \mathbf{E}[\mathbf{E}[X(I_{k,1}(\mathbf{i})) |\mathcal{M}_{ \nu_n(I_{k,1}(\mathbf{i})) \cup \nu_n(I_{k,2}(\mathbf{i}))},X(I_{k,2}(\mathbf{i}))] X(I_{k,2}(\mathbf{i})) |\mathcal{M}_{ \nu_n(\mathbf{i})}]\\ \notag
	&=& \mathbf{E}[\mathbf{E}[ X(I_{k,1}(\mathbf{i})) |\mathcal{M}_{ \nu_n(I_{k,1}(\mathbf{i}))}] X(I_{k,2}(\mathbf{i})) |\mathcal{M}_{ \nu_n(\mathbf{i})}]\\ \notag
	&=& \mathbf{E}[ X(I_{k,1}(\mathbf{i})) |\mathcal{M}_{ \nu_n(I_{k,1}(\mathbf{i}))}] \mathbf{E}[ X(I_{k,2}(\mathbf{i})) |\mathcal{M}_{ \nu_n(\mathbf{i})}]\\ \notag
	&=& \mathbf{E}[ X(I_{k,1}(\mathbf{i})) |\mathcal{M}_{ \nu_n(I_{k,1}(\mathbf{i}))}] \mathbf{E}[ X(I_{k,2}(\mathbf{i})) |\mathcal{M}_{ \nu_n(I_{k,2}(\mathbf{i}))}].
\end{eqnarray*}
The second equality follows by (\ref{nu}). The third equality follows because $\nu_n(I_{k,2}(\mathbf{i})) \setminus \bar \nu_n(I_{k,1}(\mathbf{i}))$ and $I_{k,2}(\mathbf{i})$ are outside $\bar \nu_n(I_{k,1}(\mathbf{i}))$ and $X(I_{k,1}(\mathbf{i}))$ is conditionally independent of $(\mathcal{M}_{\nu_n(I_{k,2}(\mathbf{i})) \setminus \bar \nu_n(I_{k,1}(\mathbf{i}))},X(I_{k,2}(\mathbf{i})))$ given $\mathcal{M}_{ \nu_n(I_{k,1}(\mathbf{i}))}$ by the CND property of $\{X_i\}$. The fourth equality follows because $\nu_n(I_{k,1}(\mathbf{i})) \subset \nu_n(\mathbf{i})$ . The fifth equality uses the fact that $\nu_n(I_{k,1}(\mathbf{i})) \setminus \bar \nu_n(I_{k,2}(\mathbf{i}))$ is outside of $\bar \nu_n(I_{k,2}(\mathbf{i}))$ and the CND property of $\{X_i\}_{i \in N_n}$.

Let us turn to the second statement of the lemma and now assume that Condition A holds. We write
\begin{eqnarray*}
	\label{eq21}
	\mathbf{E}[X(\mathbf{i}) |\mathcal{G}] 
	&=& \mathbf{E}[\mathbf{E}[X(I_{k,1}(\mathbf{i})) |\mathcal{M}_{ \nu_n(I_{k,1}(\mathbf{i})) },X(I_{k,2}(\mathbf{i}))] X(I_{k,2}(\mathbf{i})) |\mathcal{G}]\\ \notag
	&=& \mathbf{E}[\mathbf{E}[ X(I_{k,1}(\mathbf{i})) |\mathcal{M}_{ \nu_n(I_{k,1}(\mathbf{i}))}] X(I_{k,2}(\mathbf{i})) |\mathcal{G}]\\ \notag
	&=& \mathbf{E}[ X(I_{k,1}(\mathbf{i})) |\mathcal{G}] \mathbf{E}[ X(I_{k,2}(\mathbf{i})) |\mathcal{G}].
\end{eqnarray*}
The second equality follows because $I_{k,2}(\mathbf{i})$ is outside $\bar \nu_n(I_{k,1}(\mathbf{i}))$ and by the CND property. The third equality follows by Condition A, i.e., $\mathcal{M}_{ \nu_n(I_{k,1}(\mathbf{i}))}$ and $X(I_{k,2}(\mathbf{i}))$ are conditionally independent given $\mathcal{G}$. $\blacksquare$
\medskip

Let us present the proof of Theorem \ref{thm: CLT}. Recall the notation in the theorem, and define $X_{i}^\sigma \equiv X_{i}/\sigma_n$, $W^\sigma \equiv \sum_{i \in N_n} X_{i}^\sigma$,  and $W_{i}^\sigma \equiv \sum_{j \in \bar \nu_n (i)} X_{j}^\sigma$. Let us define for $t \in \mathbf{R}$ and $\omega \in \Omega$,
\begin{eqnarray*}
	\varphi_n(t)(\omega) \equiv \mathbf{E}[\exp(\text{i}tW^\sigma)|\mathcal{G}](\omega),
\end{eqnarray*}
where $\text{i} \equiv \sqrt{-1}$. Note that $\varphi_n$ is uniformly continuous on $\mathbf{R}$ almost surely, and since $\mathbf{E}[(W^\sigma)^2|\mathcal{G}] = 1, \text{a.e.}$, $\varphi_n$ is twice continuously differentiable almost surely. (See \cite{Yuan/Lei:16:CSTM}.)

\begin{lemma}
	\label{lemm: characteristic}
	Suppose that the conditions of Theorem \ref{thm: CLT} hold. Then for each $t \in \mathbf{R}$,
	\begin{eqnarray*}
		|\varphi_n'(t) + t \varphi_n (t)|
		\lesssim \left(n t^2 d_{mx} d_{av} \mu_3^3 + t \sqrt{n d_{mx}^2 d_{av} \mu_4^4 + r_n^2} \right), \text{a.e.}
	\end{eqnarray*}
\end{lemma}

\noindent \textbf{Proof:} First, as in the proof of Theorem 1 in \cite{Jenish/Prucha:09:JOE}, we decompose
\begin{eqnarray*}
	(\text{i}t - W^\sigma) \exp(\text{i}tW^\sigma) = h_{1,t}(W^\sigma) - h_{2,t}(W^\sigma) - h_{3,t}(W^\sigma),
\end{eqnarray*}
where
\begin{eqnarray*}
	h_{1,t}(W^\sigma) &\equiv& \text{i}t e^{\text{i}t W^\sigma} \left(1 - \sum_{j \in N_n} X_j^\sigma W_j^\sigma\right),\\
	h_{2,t}(W^\sigma) &\equiv& e^{\text{i}t W^\sigma} \sum_{j \in N_n} X_j^\sigma(1 - \text{i}t W_j^\sigma - e^{-\text{i}t W_j^\sigma}), \textnormal{ and }\\
	h_{3,t}(W^\sigma) &\equiv& e^{\text{i}t W^\sigma} \sum_{j \in N_n} X_j^\sigma e^{-\text{i}t W_j^\sigma}.
\end{eqnarray*}
Now, let us consider 
\begin{eqnarray*}
	\mathbf{E}[h_{1,t}(W^\sigma)^2|\mathcal{G}]  &\le& t^2 \mathbf{E}\left[\left(1 - \sum_{i \in N_n} X_i^\sigma W_i^\sigma  \right)^2 |\mathcal{G} \right]\\
	&=&  
	t^2 \mathbf{E}\left[\left(\sum_{i \in N_n} \sum_{j \in \bar \nu_n (i)} X_{i}^\sigma X_{j}^\sigma \right)^2 |\mathcal{G} \right] - t^2, 
\end{eqnarray*}
because
\begin{eqnarray*}
	\sum_{i \in N_n} \sum_{j \in \bar \nu_n (i)} \mathbf{E}\left[X_{i}^\sigma X_{j}^\sigma |\mathcal{G} \right] = 1,
\end{eqnarray*}
by the definition of $\sigma_n^2$.

Define $\sigma_{ij} \equiv \mathbf{E}[X_{i}^\sigma X_{j}^\sigma |\mathcal{G}]$ and $\sigma_{ij}^* \equiv \mathbf{E}[X_{i}^\sigma X_{j}^\sigma |\mathcal{M}_{\nu_n(i,j)}]$. Note that for $\{i,j\},\{k,l\} \subset N_n$ such that $\{k,l\} \subset N_n \setminus\bar \nu_n(i,j)$ and $\{i,j\} \subset N_n \setminus\bar \nu_n(k,l)$, we have by Lemma \ref{lemm: eq4},
\begin{eqnarray}
\label{eq64}
\mathbf{E}[X_{i}^\sigma X_{j}^\sigma X_{k}^\sigma X_{l}^\sigma |\mathcal{M}_{ \nu_n(i,j,k,l)}]
&=& \mathbf{E}[X_{i}^\sigma X_{j}^\sigma |\mathcal{M}_{\nu_n(i,j)}] \mathbf{E}[  X_{k}^\sigma X_{l}^\sigma |\mathcal{M}_{\nu_n(k,l)}] \\ \notag
&=& \sigma_{ij}^* \sigma_{kl}^*.
\end{eqnarray}
Let
\begin{eqnarray*}
	\tilde \Delta_n \equiv \left(\sum_{i \in N_n} \sum_{j \in \bar \nu_n(i)} \sigma_{ij}\right)^2 - \left(\sum_{i \in N_n} \sum_{j \in \bar \nu_n(i)} \sigma_{ij}^*\right)^2.
\end{eqnarray*}
Then we can write
\begin{eqnarray}
\label{eq34}
1 &=& \left(\sum_{i \in N_n} \sum_{j \in \bar \nu_n(i)} \sigma_{ij}\right)^2 \\ \notag
&=& \left(\sum_{i \in N_n} \sum_{j \in \bar \nu_n(i)} \sigma_{ij}^*\right)^2 + \tilde \Delta_n
= \sum' \sigma_{ij}^*\sigma_{kl}^* + \sum'' \sigma_{ij}^* \sigma_{kl}^* + \tilde \Delta_n,
\end{eqnarray}
where the sum $\sum'$ is over $(i,j,k,l)$ such that $i \in N_n$, $j \in \bar \nu_n(i)$, $k \in N_n$, $l \in \bar \nu_n(k)$ and either $\{k,l\} \cap \bar \nu_n(i,j) \ne \varnothing$ or $\{i,j\} \cap \bar \nu_n(k,l) \ne \varnothing$, and the sum $\sum''$ is over $(i,j,k,l)$ such that $i \in N_n$, $j \in \bar \nu_n(i)$, $k \in N_n$, $l \in \bar \nu_n(k)$ and $\{k,l\} \subset N_n \setminus \bar \nu_n(i,j)$ and $\{i,j\} \subset N_n \setminus \bar \nu_n(k,l)$. This implies that
\begin{eqnarray*}
	&& \mathbf{E}\left[\left(\sum_{i \in N_n} \sum_{j \in \bar \nu_n (i)} X_{i}^\sigma X_{j}^\sigma \right)^2 |\mathcal{G} \right]-1\\
	&=& \sum' \mathbf{E}[X_{i}^\sigma X_{j}^\sigma X_{k}^\sigma X_{l}^\sigma |\mathcal{G}] 
	+ \sum'' \mathbf{E}[X_{i}^\sigma X_{j}^\sigma X_{k}^\sigma X_{l}^\sigma |\mathcal{G}] -1 \\
	&=& \sum' \mathbf{E}[X_{i}^\sigma X_{j}^\sigma X_{k}^\sigma X_{l}^\sigma |\mathcal{G}] 
	+ \sum'' \mathbf{E}[\sigma_{ij}^*\sigma_{kl}^*|\mathcal{G}] - 1\\
	&=&  \sum' \left(\mathbf{E}[X_{i}^\sigma X_{j}^\sigma X_{k}^\sigma X_{l}^\sigma |\mathcal{G}] -    \mathbf{E}[\sigma_{ij}^*\sigma_{kl}^*|\mathcal{G}] \right) - \mathbf{E}[\tilde \Delta_n|\mathcal{G}].
\end{eqnarray*}
The second equality is by (\ref{eq64}) and the third equality is by (\ref{eq34}). The leading sum is bounded by $Cnd_{mx}^2 d_{av} \mu_4^4$, because the number of the terms in the sum $\sum'$ is bounded by $C_1 n d_{mx}^2 d_{av}$ for some constant $C_1>0$. 

Let us focus on $\tilde \Delta_n$. We write $\tilde \Delta_n=R_n(2-R_n)$, by using (\ref{eq34}), where
\begin{eqnarray*}
	R_n &\equiv& \sum_{i \in N_n} \sum_{j \in \bar \nu_n(i)} (\sigma_{ij}-\sigma_{ij}^*).
\end{eqnarray*}
Since $\mathbf{E}[R_n|\mathcal{G}] = 0$, we have 
\begin{eqnarray*}
	\mathbf{E}[\tilde \Delta_n|\mathcal{G}] = -\mathbf{E}[R_n^2|\mathcal{G}] = - r_n^2.
\end{eqnarray*}
Hence,
\begin{eqnarray*}
	\mathbf{E}[h_{1,t}(W^\sigma)^2|\mathcal{G}] \lesssim t^2 \left(n d_{mx}^2 d_{av} \mu_4^4 + r_n^2\right).
\end{eqnarray*}

Let us turn to $h_{2,t}(W^\sigma)$. Using series expansion of $\exp(-\text{i}tx)$ (e.g. see  (3.2) of \cite{Tikhomirov:80:TPA}), we bound 
\begin{eqnarray*}
	\mathbf{E}[h_{2,t}(W^\sigma)|\mathcal{G}] \le \frac{t^2}{2}\sum_{i \in N_n} \mathbf{E}[|X_i^\sigma|(W_i^\sigma)^2|\mathcal{G}]
	= \frac{t^2}{2} \sum_{i \in N_n} \sum_{j,k \in \bar \nu_n (i)} \mathbf{E}[|X_{i}^\sigma X_{j}^\sigma X_{k}^\sigma ||\mathcal{G}].
\end{eqnarray*}
Using arithmetic-geometric mean inequality, we can bound the last term by $C t^2 nd_{mx}d_{av}\mu_{3}^3$.

Finally, let us turn to $h_{3,t}(W^\sigma)$. We write $\mathbf{E}[h_{3,t}(W^\sigma)|\mathcal{G}]$ as
\begin{eqnarray*}
	\sum_{i \in N_n}\mathbf{E}[X_i^\sigma \exp(\text{i}t(W^\sigma - W_i^\sigma)) |\mathcal{G}] 
	= \sum_{i \in N_n} \mathbf{E}\left[X_i^\sigma \exp \left(\text{i}t\sum_{j \in N_n \setminus \bar \nu_n(i)} X_j^\sigma \right)|\mathcal{G} \right].
\end{eqnarray*}
The last conditional expectation is equal to
\begin{eqnarray*}
	&& \mathbf{E}\left[ \mathbf{E}\left[X_i^\sigma \exp \left(\text{i}t\sum_{j \in N_n \setminus \bar \nu_n(i)} X_j^\sigma \right) |\mathcal{M}_{\nu_n(i)}\right] |\mathcal{G}\right]\\
	&=& \mathbf{E}\left[ \mathbf{E}\left[X_i^\sigma |\mathcal{M}_{\nu_n(i)}\right] 
	\mathbf{E}\left[ \exp \left(\text{i}t\sum_{j \in N_n \setminus \bar \nu_n(i)} X_j^\sigma \right) |\mathcal{M}_{\nu_n(i)}\right] |\mathcal{G}\right] = 0.
\end{eqnarray*}
The first equality follows by CND and the second equality follows because 
$\mathbf{E}\left[X_i^\sigma |\mathcal{M}_{\nu_n(i)}\right] = 0$. Hence, it follows that
\begin{eqnarray*}
	\mathbf{E}[h_{3,t}(W^\sigma)|\mathcal{G}] = 0.
\end{eqnarray*}
Since we have 
\begin{eqnarray*}
	\varphi_n'(t) + t \varphi_n (t) = - \text{i}\left(\mathbf{E}[(\text{i}t - W^\sigma) \exp(\text{i}tW^\sigma)|\mathcal{G}]\right),
\end{eqnarray*}
by collecting the results for $h_{1,t}(W^\sigma)$, $h_{2,t}(W^\sigma)$, and $h_{3,t}(W^\sigma)$, we obtain the desired result. $\blacksquare$
\medskip

\noindent \textbf{Proof of Theorem \ref{thm: CLT}:}  For each $t \in \mathbf{R}$,
\begin{eqnarray*}
	\varphi_n'(t) = -t\varphi_n(t) - \text{i} \mathbf{E}[(\text{i}t - W^\sigma)\exp(\text{i}tW^\sigma)|\mathcal{G}] \equiv -t\varphi_n(t) + \gamma_n(t), \textnormal{say}.
\end{eqnarray*}
Taking integral of both sides, we obtain the following expression:
\begin{eqnarray*}
	\varphi_n(t) = \exp \left(- \frac{t^2}{2}\right) \left[ 1+ \int_0^t \gamma_n(u)\exp\left(\frac{u^2}{2} \right) du \right]
\end{eqnarray*}
or 
\begin{eqnarray}
\label{bound}
	\left|\varphi_n(t) - \exp \left(- \frac{t^2}{2}\right) \right|
	\le \exp \left(- \frac{t^2}{2}\right) \int_0^t |\gamma_n(u)| \exp\left(\frac{u^2}{2} \right) du.
\end{eqnarray}
Note that for all $t \ge 0$,
\begin{eqnarray*}
	\int_0^t u^2 \exp(u^2/2)du \le \exp(t^2/2)t,
\end{eqnarray*}
and $\int_0^t u \exp(u^2/2)du = \exp(t^2/2)-1$. Applying Lemma \ref{lemm: characteristic}, the last term in (\ref{bound}) for $t>0$ is bounded by $C(t a_n + (1- \exp(-t^2/2))b_n)$, where
\begin{eqnarray*}
	a_n \equiv n d_{mx} d_{av} \mu_3^3, \text{ and } b_n \equiv \sqrt{n d_{mx}^2 d_{av}\mu_4^4 + r_n^2},
\end{eqnarray*}
for some absolute constant $C>0$. Hence for any $T \ge 1$,
\begin{eqnarray*}
    \int_{[-T,T]}\left|\frac{\varphi_n(t) - e^{-t^2/2}}{t} \right| dt
	&\le&  C a_n \int_{[-T,T]}dt
	+ C b_n \int_{[-T,T]} \left|\frac{1 - e^{-t^2/2}}{t}\right| dt.
\end{eqnarray*}
The last sum is bounded by $2 C T a_n +  2 C \log(T) b_n$. Therefore, by Esseen's inequality (see e.g. Theorem 1.5.2 of \cite{Ibragimov/Linnik:71:ISSRV}, p.27), we obtain the following bound on the event $\mathcal{A}_n$,
\begin{eqnarray*}
	\Delta_n(t;\mathcal{G}) \lesssim \left(T a_n + \log(T) b_n + T^{-1} \right)
	\lesssim (a_n^{1/2} - \log(a_n) b_n),
\end{eqnarray*}
by taking $T = a_n^{-1/2}.$
$\blacksquare$
\medskip

\noindent \textbf{Proof of Lemma \ref{CLT cond}:} For $(i,j)$ and $(i',j')$ such that either $\{i,j\} \cap \bar \nu_n(i',j') = \varnothing$ or $\{i',j'\} \cap \bar \nu_n(i,j) = \varnothing$, $\mathbf{E}[\xi_{ij}\xi_{i'j'}|\mathcal{G}] = \mathbf{E}[\xi_{ij}|\mathcal{G}]\mathbf{E}[\xi_{i'j'}|\mathcal{G}]=0$. Let $A$ be the set of $((i,j),(i',j'))$ such that $i \in N_n$, $j \in \bar \nu_n(i), i' \in N_n$ and $j' \in \bar \nu_n(i')$. Then
\begin{eqnarray*}
	\mathbf{E}\left[\left(\sum_{i \in N_n} \sum_{j \in \bar \nu_n(i)} \xi_{ij} \right)^2 |\mathcal{G}\right] \le \sigma_n^4 \mu_4^4 \sum_{j=1}^4 A_j, 
\end{eqnarray*} 
where $A_1$ is the number of $((i,j),(i',j')) \in A$ such that either $i \in \bar \nu_n(i')$ or $i' \in \bar \nu_n(i)$; $A_2$ is the number of $((i,j),(i',j')) \in A$ such that either $i \in \bar \nu_n(j')$ or $j' \in \bar \nu_n(i)$; $A_3$ is the number of $((i,j),(i',j')) \in A$ such that either $j \in \bar \nu_n(i')$ or $i' \in \bar \nu_n(j)$; $A_4$ is the number of $((i,j),(i',j')) \in A$ such that either $j \in \bar \nu_n(j')$ or $j' \in \bar \nu_n(j)$. Thus, it is not hard to see that
\begin{eqnarray*}
	\sum_{j=1}^4 A_j \le 8 n d_{mx}^2 d_{av},
\end{eqnarray*}
completing the proof. $\blacksquare$
\medskip

\noindent \textbf{Proof of Corollary \ref{CLT2}:} Similarly as in the proof of Lemma \ref{lemm: characteristic}, we decompose
\begin{eqnarray*}
		(\text{i}t - W^\sigma) \exp(\text{i}tW^\sigma) = h_{1,t}(W^\sigma) - h_{2,t}(W^\sigma) - h_{3,t}(W^\sigma).
\end{eqnarray*}
The treatment of $h_{2,t}(W^\sigma)$ and $h_{3,t}(W^\sigma)$ is the same as that of the proof of Lemma \ref{lemm: characteristic}. The difference lies in the treatment of $h_{1,t}(W^\sigma)$. Using Condition A and Lemma \ref{lemm: eq4}, we note that for $\{i,j\},\{k,l\} \subset N_n$ such that $\{k,l\} \subset N_n \setminus\bar \nu_n(i,j)$ and $\{i,j\} \subset N_n \setminus\bar \nu_n(k,l)$,
\begin{eqnarray}
\label{eq61}
\mathbf{E}[X_{i}^\sigma X_{j}^\sigma X_{k}^\sigma X_{l}^\sigma |\mathcal{G}] 
= \mathbf{E}[X_{i}^\sigma X_{j}^\sigma |\mathcal{G}] \mathbf{E}[  X_{k}^\sigma X_{l}^\sigma |\mathcal{G}] = \sigma_{ij} \sigma_{kl}.
\end{eqnarray}
Following the same argument in the proof of Lemma \ref{lemm: characteristic}, we find that
\begin{eqnarray*}
	\mathbf{E}\left[\left(\sum_{i \in N_n} \sum_{j \in \bar \nu_n (i)} X_{i}^\sigma X_{j}^\sigma \right)^2 |\mathcal{G} \right]-1
	=  \sum' \left(\mathbf{E}[X_{i}^\sigma X_{j}^\sigma X_{k}^\sigma X_{l}^\sigma |\mathcal{G}] -    \sigma_{ij}\sigma_{kl}\right),
\end{eqnarray*}
which is bounded by $Cnd_{mx}^2 d_{av} \mu_4^4$. Hence in the proof of Lemma \ref{lemm: characteristic}, we do not need to deal with $\tilde \Delta_n$. Following the proofs of Lemma \ref{lemm: characteristic} and Theorem \ref{thm: CLT} for the rest of the terms, we obtain the desired result. $\blacksquare$
\medskip

\begin{lemma}
	\label{lemm: approx bd}
	Suppose that $X$ and $Y$ are random variables such that $\mathbf{E}[Y] = 0$ and $(\mathbf{E}[|Y|^r])^{1/r} \le M$ for some constants $r,M>0$, and $F$ is the CDF on $\mathbf{R}$ with density function $f$. Then for any $t \in \mathbf{R}$,
	\begin{eqnarray*}
		\left|P\{ X + Y \le t\} - F(t) \right|
		&\le& \left|P\{ X \le t+ q\} - F(t+ q) \right|\\
		&& +\left|P\{ X \le t- q\} - F(t- q) \right|  
	     + 4 (\overline c M)^{r/(1+r)},
	\end{eqnarray*}
	where $\overline c \equiv \sup_{z \in \mathbf{R}}f(z) \text{ and }
		q \equiv (M^r/\overline c)^{1/(1+r)}$.
\end{lemma}

\noindent \textbf{Proof: } First, note that for any $\varepsilon > 0$,
\begin{eqnarray}
    \label{ineq4}
	\left|P\{X+Y \le t\} - F(t) \right| \le \left|P\{X+Y \le t, |Y| \le \varepsilon \} - F(t) \right| + P\{|Y| > \varepsilon\}.
\end{eqnarray}
As for the probability inside the absolute value above, we note that
\begin{eqnarray*}
	P\{X+Y \le t, |Y| \le \varepsilon \} &\ge& P\{X +\varepsilon \le t, |Y| \le \varepsilon \}\\
	 &\ge& P\{X +\varepsilon \le t\} - P\{|Y| > \varepsilon\}.
\end{eqnarray*}
Also, observe that
\begin{eqnarray*}
	P\{X+Y \le t, |Y| \le \varepsilon \} \le P\{X - \varepsilon \le t \}.
\end{eqnarray*}
Hence
\begin{eqnarray*}
	&& \left|P\{X+Y \le t, |Y| \le \varepsilon \} - F(t) \right| \\
	&\le& \max \left\{|P\{X +\varepsilon \le t\} - P\{|Y| > \varepsilon\} -F(t)|, |P\{X - \varepsilon \le t \} - F(t)| \right\}\\
	&\le& |P\{X \le t-\varepsilon\} -F(t-\varepsilon)| + |P\{X \le t+\varepsilon\} -F(t+\varepsilon)|\\
	&& + 2 \sup_{z \in \mathbf{R}} f(z) \varepsilon + P\{|Y| > \varepsilon\}.
\end{eqnarray*}
From (\ref{ineq4}),
\begin{eqnarray*}
	\left|P\{X+Y \le t\} - F(t) \right| &\le&  |P\{X \le t-\varepsilon\} -F(t-\varepsilon)| + |P\{X \le t+\varepsilon\} -F(t+\varepsilon)|\\
	&& + 2 \sup_{z \in \mathbf{R}} f(z) \varepsilon + 2 P\{|Y| > \varepsilon\}.
\end{eqnarray*}
Using Markov's inequality, we bound the last term by $2 \varepsilon^{-r} \mathbf{E}[|Y|^r] \le 2 \varepsilon^{-r} M^r$. Taking $\varepsilon = q$, we obtain the desired result. $\blacksquare$
\medskip

\noindent \textbf{Proof of Theorem \ref{thm: CLT-CHD}: } We write $\Delta^*(t;\mathcal{G})$ as
\begin{eqnarray}
    \label{prob diff}
    \left| P\left\{\frac{1}{\sigma_n^*} \sum_{i \in N_n^*} X_i + R_n^* \le t |\mathcal{G} \right\} - \Phi(t) \right|,
\end{eqnarray} 
where $R_n^* \equiv \frac{1}{\sigma_n^*} \sum_{i \in N_n \backslash N_n^*} X_i.$ We write
\begin{eqnarray*}
	R_n^* =  \sum_{i \in N_n\setminus N_n^*} \frac{X_i}{\sigma_n^*}  = \sum_{i \in N_n\setminus N_n^*} \xi_i, \text{ say}.
\end{eqnarray*}
Now, choose $ 1 \le r \le 4$ and write
\begin{eqnarray*}
	\mathbf{E}\left[ \left| \sum_{i \in N_n\setminus N_n^*} \xi_i \right|^r |\mathcal{G} \right]
	\le \left( \sum_{i \in N_n\setminus N_n^*} \left( \mathbf{E}[|\xi_i|^r|\mathcal{G}] \right)^{1/r} \right)^r
	\le \left( |N_n\setminus N_n^*| \tilde \mu_r^* \right)^r.
\end{eqnarray*}
Hence on the event $\mathcal{A}_{n,r}(\varepsilon_n)$,
\begin{eqnarray*}
	(\mathbf{E}[|R_n^*|^r|\mathcal{G}])^{1/r} + \rho_r^* \le \varepsilon_n.
\end{eqnarray*}
Define for brevity,
\begin{eqnarray*}
	S_{n,\sigma}^* \equiv \frac{1}{\sigma_n^*} \sum_{i \in N_n^*} (X_i - \mathbf{E}[X_i|\mathcal{M}_{\nu_n^*(i)}^*]).
\end{eqnarray*}
By Lemma \ref{lemm: approx bd}, the term (\ref{prob diff}) is bounded by (for any $r>0$)
\begin{eqnarray}
    \label{bd1}
	&& \left| P\left\{ S_{n,\sigma}^* \le t + q_n|\mathcal{G}\right\} - \Phi(t + q_n) \right|\\ \notag
	&& +\left| P\left\{ S_{n,\sigma}^* \le t - q_n|\mathcal{G}\right\} - \Phi(t - q_n) \right|
	+ \quad 4 \left( \phi(0) \varepsilon_n\right)^{r/(r+1)},
\end{eqnarray}
where $\phi$ denotes the density of $N(0,1)$ and
\begin{eqnarray*}
	q_n = \varepsilon_n^{r/(r+1)}/\phi(0)^{1/(r+1)}.
\end{eqnarray*}
By Lemma \ref{lemm: invar}, $\{X_i\}_{i \in N_n^*}$ is CND with respect to $(\nu_n^*,\mathcal{M}^*)$, we apply Theorem \ref{thm: CLT} to the leading two terms in (\ref{bd1}) to obtain their bound as
\begin{eqnarray*}
	C \sqrt{n^* d_{mx}^* d_{av}^* \mu_3^{*3}} - C \log(n^* d_{mx}^* d_{av}^*\mu_3^{*3}) \sqrt{n^* d_{mx}^{*2} d_{av}^*\mu_4^{*4} + r_n^{*2}}
\end{eqnarray*}
for some constant $C>0$, delivering the desired result. $\blacksquare$

\begin{lemma}
\label{lemm: aux}	
	Suppose that $(\mathbb{D},d)$ is a given metric space and for each $n \ge 1$, $\xi_n,\zeta_n,\zeta$ are $\mathbb{D}$-valued random variables. If for each $\varepsilon>0$, $P^*\{d(\xi_n,\zeta_n) > \varepsilon\} \rightarrow 0$ and $\zeta_n \rightarrow \zeta$, $\mathcal{G}$-stably, as $n \rightarrow \infty$, then
	\begin{eqnarray*}
		\xi_n \rightarrow \zeta, \mathcal{G}\text{-stably.}
	\end{eqnarray*} 
\end{lemma}

\noindent \textbf{Proof: } 	First note that $\zeta_n \rightarrow \zeta$, $\mathcal{G}$-stably if and only if for all event $U \in \mathcal{G}$ and any closed set $F \in \mathcal{B}(\mathbb{D})$,
\begin{eqnarray}
\label{Portmanteau}
\text{limsup}_{n \rightarrow \infty} P^*\{\zeta_n \in F\} \cap U \le P\{\zeta \in F\} \cap U.
\end{eqnarray}
(This can be shown following the proof of Theorem 1.3.4 (iii) of \cite{vanderVaart/Wellner:96:WeakConvg}.) Using this and following the same arguments in the proof of Lemma 1.10.2 of \cite{vanderVaart/Wellner:96:WeakConvg}, we deduce the lemma. $\blacksquare$
\medskip

\noindent \textbf{Proof of Lemma \ref{generalized slutsky}: } (i) Since $P^*\{d(\xi_n,\xi)>\varepsilon\} \rightarrow 0$, we have $P^*\{\tilde d((\zeta_n,\xi_n),(\zeta_n,\xi)) >\varepsilon\} \rightarrow 0$, where $\tilde d$ is a metric on $\mathbb{D} \times \mathbb{D}$ defined as $\tilde d((f_1,f_2),(g_1,g_2)) = d(f_1,g_1) + d(f_2,g_2)$ for $f_1,f_2,g_1,g_2 \in \mathbb{D}$. Furthermore, note that $(\zeta_n,\xi) \rightarrow (\zeta,\xi)$, $\mathcal{G}$-stably, because $\zeta_n \rightarrow \zeta$, $\mathcal{G}$-stably, and $\xi$ is $\mathcal{G}$-measurable. Now the desired result follows by Lemma \ref{lemm: aux}.
\medskip

\noindent (ii) Note that $\zeta_n' \rightarrow \zeta'$, $\mathcal{G}$-stably, if and only if for any event $U \in \mathcal{G}$ and any open set $G \in \mathcal{B}(\mathbb{D})$,
\begin{eqnarray}
\label{Portmanteau2}
\text{liminf}_{n \rightarrow \infty} P_*\{\zeta_n' \in G\} \cap U \ge P\{\zeta' \in G\} \cap U,
\end{eqnarray}
where $P_*$ denotes the inner probability. Using this and following the same arguments in the proof of Theorem 3.27 of \cite{Kallenberg:1997:Foundations} for the continuous mapping theorem for weak convergence, we obtain the proof of (ii). $\blacksquare$
\medskip

For the proof of Theorem \ref{thm: fundamental theorem}, we use the following lemma. 
\begin{lemma}
	\label{aux lemma}
	If $f:\mathbb{D} \mapsto \mathbf{R}$ is
	bounded and continuous, and $K\subset \mathbb{D}$ is compact, then for every $\epsilon >0$ there exists $\tau >0$
	such that, if $x\in K$ and $y\in \mathbb{D}$
	with $\left\Vert x-y\right\Vert <\tau $, then 
	\begin{equation*}
		\left\vert f(x)-f(y)\right\vert <\epsilon \text{.}
	\end{equation*}
\end{lemma}

\noindent \textbf{Proof of Theorem \ref{thm: fundamental theorem}:} First, let us suppose that (i) and (ii) hold. To see that the marginal $PL$ of $L$ is a tight Borel law, note that the stable finite dimensional convergence of $\zeta_n$ implies the convergence of the finite dimensional distributions of $\zeta_n$. Combining this with the asymptotic $\rho$-equicontinuity and using Theorems 1.5.4 and 1.5.7 of \cite{vanderVaart/Wellner:96:WeakConvg}, we obtain that $PL$ is a tight Borel law. The fact that $PL$ is concentrated on $U_\rho(\mathcal{H})$ follows from Theorem 10.2 of \cite{Pollard:90:EmpiricalProcesses}. 

Now let us show the $\mathcal{G}$-stable convergence of $\zeta_n$. We follow the arguments in the proof of Theorem 2.1 of \cite{Wellner:05:EmpTheory}. Let $\zeta \in U_\rho(\mathcal{H})$ be a random element whose distribution is the same as $PL$. Since $\left(\mathcal{H},\rho \right)$ is
totally bounded, for every $\delta >0$ there exists a finite set of points 
$\mathcal{H}_\delta$ that is $\delta $-dense in $\mathcal{H}$ i.e. $\mathcal{H\subset \cup }%
_{h \in \mathcal{H}_\delta} B\left(h,\delta \right) $ where $
B\left(h,\delta \right) $ is the open ball with center $h$ and radius $
\delta $. Thus, for each $h\in \mathcal{H}$, we can choose $\pi _{\delta
}\left(h\right) \in \mathcal{H}_\delta$
such that $\rho( \pi_{\delta }(h),h) <\delta $. Define
\begin{equation*}
	\zeta_{n,\delta }\left(h\right) = \zeta_{n}\left( \pi_{\delta }(h) \right) \text{, and } \zeta_{\delta }(h)= \zeta(\pi_\delta(h)) \text{ for } h \in \mathcal{H}.
\end{equation*}%
By the $\mathcal{G}$-stable
convergence of the finite dimensional projection of $\zeta _{n}$, we have for
each $U \in \mathcal{G}$ and for each
bounded and continuous functional $f:\ell ^{\infty }\left( \mathcal{H}%
\right) \mapsto \mathbf{R}$, 
\begin{eqnarray}
\label{conv1}
	\mathbf{E}^*\left[ f(\zeta_{n,\delta }) 1_U \right] \rightarrow
	\mathbf{E}\left[f(\zeta_\delta )1_U \right].
\end{eqnarray}
Furthermore, the a.e. uniform continuity of the sample paths of $\zeta $ implies that 
\begin{equation}
\label{conv2}
P\left\{	\lim_{\delta \rightarrow 0}\sup_{h \in \mathcal{H}} \left| \zeta(h) -\zeta _{\delta }(h)\right| = 0 \right\} = 1.
\end{equation}%
For each bounded and continuous functional $f:\ell ^{\infty }\left( \mathcal{H}\right) \mapsto \mathbf{R}$, and for each $U \in \mathcal{G}$ such that $P(U)>0$, 
\begin{eqnarray*}
	&& |\mathbf{E}^*\left[ f(\zeta _{n})1_U\right] -\mathbf{E}\left[
	f(\zeta )1_U\right] | \\
	&\leq & |\mathbf{E}^*\left[ (f(\zeta _{n})-f(\zeta _{n,\delta }))1_U\right]| + | \mathbf{E}\left[ (f(\zeta _{n,\delta })-f(\zeta _{\delta }))1_U\right]| + | \mathbf{E}\left[ (f(\zeta _{\delta }) -
	 f(\zeta ))1_U\right]|.
\end{eqnarray*}
The last two absolute values vanish as $n \rightarrow \infty$ and then $\delta \rightarrow 0$ by (\ref{conv1}) and by (\ref{conv2}) combined with the Dominated Convergence Theorem. We  use the asymptotic $\rho$-equicontinuity of $\zeta_n$ and Lemma \ref{aux lemma} and the fact that $PL$ is a tight law, and follow standard arguments to show that the leading difference vanishes as $n \rightarrow \infty$ and then $\delta \rightarrow 0$. (See the proof of Theorem 2.1 of \cite{Wellner:05:EmpTheory} for details.) 

Since the $\mathcal{G}$-stable convergence of $\zeta_n$ implies that of its finite dimensional distributions, and the weak convergence of $\zeta_n$, the converse can be shown using the standard arguments. (Again, see the proof of  Theorem 2.1 of \cite{Wellner:05:EmpTheory}.) $\blacksquare$
\medskip
	
\noindent \textbf{Proof of Lemma \ref{lemma: tail bound}:} The proof follows that of Theorems 2.3 and 3.4 of \cite{Janson:04:RSA}. In particular (\ref{bd3}) follows from Theorem 2.3. However, for (\ref{bd2}), we need to modify the proof of Theorem 3.4 because $\sigma_i^2$'s are not necessarily $\mathcal{G}$-measurable, and hence the equations (3.9) and (3.10) on page 241 do not necessarily follow.

First, without loss of generality, we set $M=1$. Following the proof of Theorem 3.4 in \cite{Janson:04:RSA} (see (3.7) there), we obtain that for any $c \ge 0$,
\begin{eqnarray}
\label{bd7}
	\mathbf{E}[\exp(cX_i)|\mathcal{M}_{\nu_n(i)}] \le
	\exp(\sigma_i^2 g(c)), 
\end{eqnarray}
where $g(c) = e^c - 1 - c$. Let $N_j^* \subset N_n$, $j=1,...,J$ be disjoint subsets which partition $N_n$ such that for any $i_1,i_2 \in N_j^*$, $i_1 \ne i_2$, $i_1 \notin \nu_n(i_2)$ and $i_2 \notin \nu_n(i_1)$. Fix $u \ge 0$, $p_j \ge 0, j=1,...,J$ such that $\sum_{j=1}^J p_j =1$ and $w_j \in [0,1], j=1,...,J$ such that $\sum_{j \in N_n: i \in N_j^*} w_j = 1$ for all $i \in N_n$. Then using Lemma \ref{lemm: eq4} and (\ref{bd7}) and following the same argument in (3.8) of \cite{Janson:04:RSA},
\begin{eqnarray*}
	\mathbf{E}\left[\exp\left(u\sum_{i \in N_n} X_i\right)|\mathcal{G}\right]
	&\le& \sum_{j=1}^J p_j \mathbf{E}\left[\prod_{i \in N_j^*} \mathbf{E}\left[\exp\left(\frac{w_j u}{p_j} X_i \right) |\mathcal{M}_{\nu_n(i)}\right]|\mathcal{G}\right]\\
	&\le& \sum_{j=1}^J p_j \mathbf{E}\left[\exp \left(\sum_{i \in N_j^*} \sigma_i^2 g\left(\frac{w_ju}{p_j}\right) \right)|\mathcal{G}\right]\\
	&\le& \sum_{j=1}^J p_j \mathbf{E}\left[\exp \left(\sum_{i \in N_n} \sigma_i^2 g\left(\frac{w_ju}{p_j}\right) \right)|\mathcal{G}\right],
\end{eqnarray*}
because $g(\cdot) \ge 0$. The last term above is bounded by
\begin{eqnarray}
    \label{bd77}
	\sum_{j=1}^J p_j \mathbf{E}\left[\exp \left( (\kappa_n + V_n)g\left(\frac{w_ju}{p_j}\right) \right)|\mathcal{G}\right]
\end{eqnarray}
where
\begin{eqnarray*}
	\kappa_n = \log \mathbf{E}\left[ \exp\left(\left|\sum_{i \in N_n} \sigma_i^2 - V_n\right|\right)|\mathcal{G}\right].
\end{eqnarray*}
As for $\kappa_n$, note that
\begin{eqnarray}
\label{bd51}
	\kappa_n \le \log \mathbf{E}\left[ \exp\left(\left|V_n\right|\right)|\mathcal{G}\right] = V_n,
\end{eqnarray}
because $V_n$ is $\mathcal{G}$-measurable. Let $W = \sum_{j=1}^J w_j$ and take $p_j = w_j/W$ to rewrite (\ref{bd77}) as
\begin{eqnarray*}
	\mathbf{E}\left[\exp \left((\kappa_n + V_n) g\left(uW\right) \right)|\mathcal{G}\right].
\end{eqnarray*}
Hence we have for each $t \ge 0$,
\begin{eqnarray*}
	P\left\{\sum_{i \in N_n} X_i > t |\mathcal{G}\right\}
	\le \mathbf{E}\left[\exp \left((\kappa_n + V_n) g\left(uW\right) - ut \right)|\mathcal{G}\right].
\end{eqnarray*}
If we take
\begin{eqnarray*}
	u = \frac{1}{W}\log\left(\frac{t}{(\kappa_n + V_n)W} + 1\right)
\end{eqnarray*}
and let $\varphi(x) \equiv (1+x)\log(1+x) - x$, the last bound becomes
\begin{eqnarray*}
	\exp\left(- (\kappa_n + V_n)\varphi\left( \frac{t}{(\kappa_n + V_n)W} \right) \right)
	&\le& \exp \left(- \frac{t^2}{2W(W(\kappa_n + V_n) + t/3)}\right)\\
	&\le& \exp \left(- \frac{t^2}{2W(2 W V_n + t/3)}\right),
\end{eqnarray*}
where the first inequality follows by the inequality: $\varphi(x) \ge x^2/(2(1+x/3))$, $x \ge 0$, and the last inequality follows by (\ref{bd51}). Now, as in the proof of Theorem 2.3 of \cite{Janson:04:RSA}, the rest of the proof can be proceeded by taking $\{(N_j^*,w_j)\}_{j=1}^J$ as a minimal fractional proper cover of $N_n$. $\blacksquare$
\medskip

\noindent \textbf{Proof of Lemma \ref{lemm: maximal inequality}}: We adapt the proof of Theorem A.2 of \cite{vanderVaart:96:AS} to accommodate the CND property of $\{Y_{i}\}_{i\in N_{n}}$. Fix $q_{0}$ so that
\begin{eqnarray*}
	2^{-q_{0}}\leq \bar{\rho}_{n}(H)\leq 2^{-q_{0}+1}
\end{eqnarray*}
and for each $q\ge q_0$, construct a nested sequence of
partitions $\mathcal{H=\cup }_{i=1}^{N_{q}}\mathcal{H}_{q_{i}}$ such that 
\begin{equation}
\frac{1}{n}\sum_{j\in N_{n}} \mathbf{E}^*\left(\sup_{h,g\in \mathcal{H}%
	_{q_{i}}}\left\vert h\left( Y_{j}\right) -g\left( Y_{j}\right) \right\vert
^{2}\right) <2^{-2q}\text{ for every }i=1,...,N_q\text{.}  \label{vdv1}
\end{equation}
By the definition of $\bar{\rho}_{n}(h)$ and the bracketing entropy, $N_{q}$ can be taken to satisfy
\begin{equation*}
\log N_{q}\leq \sum_{r=q_{0}}^{q}\log \left( 1+N_{[]}(2^{-r},\mathcal{H},%
\bar{\rho}_{n})\right) .
\end{equation*}%
Choose for each $q$ a fixed element $h_{q_{i}}$ from each $\mathcal{H}%
_{q_{i}}$ and set $\pi _{q}h=h_{q_{i}}$ and $\triangle _{q}h=( \sup_{h,g\in \mathcal{H}_{q_{i}}}|
h-g|)^*$, whenever $h\in \mathcal{H}_{q_{i}}$, where $(h)^*$ defines the minimal measurable cover of $h$ (\cite{Dudley:85:PBS}.) Then $(\bar{%
	\rho}_{n}\left( \triangle _{q}h\right))^2 <2^{-2q}$ from (\ref{vdv1}), and $\pi _{q}h$ and $\triangle _{q}h$ run through a set of $N_{q}$
functions as $h$ runs through $\mathcal{H}$. Define for each fixed $n$ and $%
q\geq q_0$, the following numbers and indicator functions: 
\begin{eqnarray*}
	\alpha _{q} &=&2^{-q}/\sqrt{\log N_{q+1}}, \\
	A_{q-1}h &=&1\left\{ \triangle _{q_{0}}h\leq \sqrt{n}\alpha
	_{q_{0}},...,\triangle _{q-1}h\leq \sqrt{n}\alpha _{q-1}\right\} , \\
	B_{q}h &=&1\left\{ \triangle _{q_{0}}h\leq \sqrt{n}\alpha
	_{q_{0}},...,\triangle _{q-1}h\leq \sqrt{n}\alpha _{q-1},\triangle _{q}h>%
	\sqrt{n}\alpha _{q}\right\} , \\
	B_{q_{0}}h &=&1\left\{ \triangle _{q_{0}}h>\sqrt{n}\alpha
	_{q_{0}}\right\} \text{.}
\end{eqnarray*}%
Because the partitions are nested, $A_{q}h$ and $B_{q}h$ are constant in $h$
on each of the partitioning sets $\mathcal{H}_{q_{i}}$ at level $q$. Now
decompose $h=\left( h-\pi _{q_{0}}h\right) +\pi _{q_{0}}h \equiv I+II$, say, with $I \equiv h-\pi _{q_{0}}h=(h-\pi _{q_{0}}h)$.  Then we can write 
\begin{eqnarray*}
	I &=&\left( h-\pi _{q_{0}}h\right)
	B_{q_{0}}h+\sum_{q=q_{0}+1}^{\infty }\left( h-\pi _{q}h\right)
	B_{q}h+\sum_{q=q_{0}+1}^{\infty }\left( \pi _{q}h-\pi _{q-1}h\right) A_{q-1}h
	\\
	&\equiv& I_{a}+I_{b}+I_{c}, \text{ say }.
\end{eqnarray*}
We analyze the empirical process at each of $I_{a}$, $I_{b}$ and 
$I_{c}$. 
\smallskip

\noindent \textbf{Control of }$I_{a}$\textbf{:} Let us bound $\mathbf{E}^* [\sup_{h\in \mathcal{H}}|\mathbb{G}_{n}(( h-\pi _{q_{0}}h)B_{q_{0}}h)|]$ by
\begin{eqnarray*}
	&&\mathbf{E}^*\left( \sup_{h\in \mathcal{H}}\frac{1}{\sqrt{n}}%
	\sum_{i \in N_n}(|h -\pi _{q_{0}}h|B_{q_{0}}h)(X_i) \right.  \\
	&&\left. +\sup_{h\in \mathcal{H}}\frac{1}{\sqrt{n}}\sum_{i \in N_n} \mathbf{E}%
	[(|h -\pi _{q_{0}}h|B_{q_{0}}h)(X_i)|\mathcal{M}_{\nu (i)}]\right).
\end{eqnarray*}%
Since $|h-\pi _{q_{0}}h| B_{q_{0}}h\le 2H 1\{ 2H>\sqrt{n}\alpha _{q_{0}}\} \le 4H^2/(\sqrt{n}
	\alpha _{q_{0}})$, we bound the last expression by 
\begin{eqnarray*}
	&& 8 \alpha _{q_{0}}^{-1}\bar{\rho}_{n}\left( H\right) ^{2}
	\le 8 \left( 2^{-q_{0}}/\sqrt{\log N_{q_{0}+1}}\right)
	^{-1}2^{-2q_{0}+2}\le 32 \cdot 2^{-q_{0}}\sqrt{\log N_{q_{0}+1}},
\end{eqnarray*}%
due to our choice of $q_0$ satisfying that $\bar{\rho}_{n}(H)^{2}\leq 2^{-2q_{0}+2}$.
\smallskip

\noindent \textbf{Control of }$I_{c}$\textbf{: }For $I_{c}=\sum_{q=q_{0}+1}^{\infty
}\left( \pi _{q}h-\pi _{q-1}h\right) A_{q-1}h$, there are at most $N_{q}-1$
functions $\pi _{q}h-\pi _{q-1}h$ and at most $N_{q-1}-1$ functions $A_{q-1}h
$. Since the partitions are nested, the function $\left\vert \pi _{q}h-\pi
_{q-1}h\right\vert A_{q-1}h$ is bounded by $\triangle _{q-1}hA_{q-1}h\leq 
\sqrt{n}\alpha _{q-1}$. Applying Corollary \ref{maximal inequality0} (with $J=\sqrt{n}\alpha _{q-1}
$ and $m=N_{q}-1$) to $\mathbb{G}_n ( \pi _{q}h-\pi _{q-1}h) A_{q-1}h$, 
\begin{eqnarray*}
&& \mathbf{E}\left[ \max_{h \in \mathcal{H}}|\mathbb{G}_{n}\left( \pi _{q}h-\pi _{q-1}h\right)
A_{q-1}h||\mathcal{G}\right] \\ \notag 
&& \lesssim (d_{mx}+1)\left( \alpha _{q-1}\log N_{q}+%
\sqrt{(\log N_{q}) \max_{1\leq s\leq N_{q}-1} V_{n}(h_{s})}\right) 
\end{eqnarray*}%
where $\{h_s: s=1,...,N_{q-1}\} \equiv \{ (\pi_q - \pi_{q-1})h A_{q-1} h: h \in \mathcal{H}\}$. From the law of
the iterated conditional expectations and Jensen's inequality,
\begin{eqnarray*}
	\mathbf{E}^* \left[\sup_{h \in \mathcal{H}} |\mathbb{G}_{n}I_{c}| \right]
	&\lesssim &(d_{mx}+1)\sum_{q=q_{0}+1}^{\infty }\left( \sqrt{(\log
		N_{q})\mathbf{E}^* \left[\max_{1\leq s\leq N_{q}-1} V_{n}(h_{s})\right]}+2^{-q}\sqrt{\log N_{q}}\right)  \\
	&\lesssim &2(d_{mx}+1)\sum_{q=q_{0}+1}^{\infty }2^{-q}\sqrt{\log N_{q}},
\end{eqnarray*}%
where for the last inequality, we used (\ref{vdv1}) so that  
\begin{equation}
\mathbf{E}^* \left[\max_{1\leq s\leq N_{q}-1}
V_{n}(h_{s})\right] \leq \mathbf{E}^*\left( \max_{1\leq s\leq
	N_{q}-1}\left( \frac{1}{n} \sum_{i\in N_{n}}\mathbf{E}\left[ \mathbf{E}%
[h_{s}^{2}\left( Y_{i}\right) |\mathcal{M}_{i}]|\mathcal{G}\right] \right)
\right) \lesssim 2^{-q+1}.  \label{last bound}
\end{equation}

\noindent \textbf{Control of }$I_{b}$ and $II$\textbf{: } The proof of these parts are the same
as that of Theorem A.2 of van der Vaart (1996) except that we use $\bar{\rho}%
_{n}(\cdot )$ instead of $\left\Vert \cdot \right\Vert _{P,2}$ so that we have $%
\mathbf{E}^*\left\Vert \mathbb{G}_{n}I_{b}\right\Vert _{\mathcal{H}%
}\lesssim \sum_{q=q_{0}+1}^{\infty }2^{-q}\sqrt{\log N_{q}}$ and $\mathbf{E}%
^*\left\Vert \mathbb{G}_{n}II\right\Vert _{\mathcal{H}}\lesssim
2^{-q_{0}}\sqrt{\log N_{q_{0}}}$. 
\smallskip

Now collecting the results for $I_{a}$, $I_{b}$, $I_{c}$ and $II$, we have 
\begin{eqnarray*}
	\mathbf{E}^* \left[\sup_{h \in \mathcal{H}}|\mathbb{G}_{n}(h)|\right] &\lesssim& (d_{mx}+1)\sum_{q=q_{0}+1}^{\infty }2^{-q}\sqrt{\log N_{q}} \\
	&\lesssim& (d_{mx}+1) \int_{0}^{\bar{\rho}_{n}(H)}\sqrt{1+\log
		N_{[]}(\varepsilon ,\mathcal{H},\bar{\rho}_{n})}d\varepsilon 
	\text{,}
\end{eqnarray*}%
giving the required result. $\blacksquare$
\medskip

\noindent \textbf{Proof of Theorem \ref{thm: empirical CLT}:} We prove conditions for Theorem \ref{thm: fundamental theorem}. Let us first consider the convergence of finite dimensional distributions. Without loss of generality, we consider the CLT for $\mathbb{G}_n(h)$ for some $h \in \mathcal{H}$ such that $Var(\mathbb{G}_n(h)|\mathcal{G}) >0$, a.e.. Assumption \ref{assump: conds emp CLT} (a) together with the moment condition for the envelope $H$ implies that $d_{mx} <C$ for all $n \ge 1$ for some $C>0$. We apply Theorem \ref{thm: CLT} to obtain the convergence of finite dimensional distributions. By the CND property of $(Y_i)_{i \in N_n}$ and Assumption \ref{assump: conds emp CLT}(b), Condition (i) in Theorem 3.3 is satisfied.

Let us prove asymptotic $\rho$-equicontinuity (with $\rho(h,h) \equiv \rho(h) \equiv \lim_{n \rightarrow \infty} \bar \rho_n(h)$.) Define $\mathcal{H}_{n,\delta} \equiv \{h-g: \bar \rho_n(h-g) \le \delta, h,g \in \mathcal{H} \}$. Then, by Lemma \ref{lemm: maximal inequality},
\begin{eqnarray*}
	\mathbf{E}^* \left[\sup_{h \in \mathcal{H}_{n,\delta}}\left|\mathbb{G}_n(h) \right| \right]
	\lesssim (1+d_{mx})\int_{0}^{\delta}\sqrt{1+\log N_{[]}(\varepsilon ,%
		\mathcal{H}_{n,\delta},\bar{\rho}_{n})}d\varepsilon.
\end{eqnarray*}
By noting that $\mathcal{H}_{n,\delta}$ is contained in $\mathcal{H} - \mathcal{H}$ and by Assumption \ref{assump: conds emp CLT}(a), the last bound vanishes as $\delta \rightarrow 0$ for each $n \ge 1$. Thus the asymptotic $\rho$-equicontinuity follows, by the condition for $\rho(h)$. $\blacksquare$
\medskip

\noindent \textbf{Proof of Theorem \ref{thm: empirical CLT 2}:} By Assumption \ref{assump: conds emp CLT 2}, we have
\begin{eqnarray*}
	\frac{n-n^*}{\sqrt{n}} \max_{i \in N_n} \sqrt{\mathbf{E}[H^2(Y_i)]} + \mathbf{E}^*\left[\sup_{h \in \mathcal{H}}\rho_n^*(h) \right] \rightarrow 0,
\end{eqnarray*}
as $n \rightarrow \infty$. Hence the desired result follows applying Theorem \ref{thm: empirical CLT} to $\mathbb{G}_n^*$ in (\ref{decomp}).$\blacksquare$

\bibliographystyle{econometrica}
\bibliography{Local_Dependence_A7}
\end{document}